\documentclass{amsart}
\usepackage[T1]{fontenc}
\usepackage[utf8]{inputenc}
\usepackage[english]{babel}
\usepackage{times}
\usepackage{afterpage}
\usepackage{booktabs}
\usepackage{tabularx}
\usepackage{graphicx}
\graphicspath{{img/}}
\usepackage{sidecap}
\usepackage{wrapfig}
\usepackage{subfigure}
\usepackage{setspace}
\usepackage{multicol}
\usepackage{multirow}
\usepackage{verbatim}
\usepackage{xcolor}
\usepackage{textcomp}
\usepackage{algorithm}
\usepackage{algorithmicx}
\usepackage{algpseudocode}
\usepackage{float}
\usepackage{dirtytalk}
\usepackage{cancel}
\usepackage{chngpage}

\usepackage{xurl}
\algdef{SE}[DOWHILE]{Do}{doWhile}{\algorithmicdo}[1]{\algorithmicwhile\ #1}

\usepackage{amsmath}
\usepackage{amssymb}
\usepackage{amsthm}
\usepackage{mathtools}

\usepackage[super]{cite}
\usepackage{url}
\usepackage[hidelinks]{hyperref}
\usepackage{cleveref}
\mathtoolsset{showonlyrefs}
\newcommand{\numberset}{\mathbb}
\newcommand{\N}{\numberset{N}}
\newcommand{\R}{\numberset{R}}

\newcommand{\E}{\numberset{E}}

\renewcommand{\epsilon}{\varepsilon}
\DeclarePairedDelimiter{\abs}{\lvert}{\rvert}
\DeclarePairedDelimiter{\norm}{\lVert}{\rVert}

\newtheorem{thm}{Theorem}

\newtheorem{cor}[thm]{Corollary}

\renewcommand{\arraystretch}{1.15}

\begin{document}

\title{Improvements on uncertainty quantification\\ with variational autoencoders}

\author{Andrea Tonini$^\mathrm{1,*}$, Tan Bui-Thanh$^\mathrm{2,3}$, Francesco Regazzoni$^\mathrm{1}$, Luca Dede'$^\mathrm{1}$ \and Alfio Quarteroni$^\mathrm{4,5}$}
\thanks{$^\mathrm{1}$MOX, Dipartimento di Matematica, Politecnico di Milano, Milan, Italy.}
\thanks{$^\mathrm{2}$Department of Aerospace Engineering and Engineering Mechanics, the University of Texas at Austin, Austin, Texas, USA.}
\thanks{$^\mathrm{3}$The Oden Institute for Computational Engineering and Sciences, the University of Texas at Austin, Austin, Texas, USA.}
\thanks{$^\mathrm{4}$(Professor Emeritus) MOX, Dipartimento di Matematica, Politecnico di Milano, \\
Piazza Leonardo da Vinci, Milan, 20133, Italy}
\thanks{$^\mathrm{5}$(Professor Emeritus) Institute of Mathematics, Ecole Polytechnique Fédérale de Lausanne,\\
 Lausanne, Switzerland}

\thanks{$^\mathrm{*}$Corresponding author. E-mail: andrea.tonini@polimi.it}

\markboth{Andrea Tonini et al.}{Uncertainty quantification variational autoencoder}

\begin{abstract}
Inverse problems aim to determine model parameters of a mathematical problem from given observational data. Neural networks can provide an efficient tool to solve these problems. In the context of Bayesian inverse problems, Uncertainty Quantification Variational AutoEncoders (UQ-VAE), a class of neural networks, approximate the posterior distribution mean and covariance of model parameters. This allows for both the estimation of the parameters and their uncertainty in relation to the observational data. In this work, we propose a novel loss function for training UQ-VAEs, which includes, among other modifications, the removal of a sample mean term from an already existing one. This modification improves the accuracy of UQ-VAEs, as the original theoretical result relies on the convergence of the sample mean to the expected value (a condition that, in high dimensional parameter spaces, requires a prohibitively large number of samples due to the curse of dimensionality). Avoiding the computation of the sample mean significantly reduces the training time in high dimensional parameter spaces compared to previous literature results. Under this new formulation, we establish a new theoretical result for the approximation of the posterior mean and covariance for general mathematical problems. We validate the effectiveness of UQ-VAEs through three benchmark numerical tests: a Poisson inverse problem, a non affine inverse problem and a 0D cardiocirculatory model, under the two clinical scenarios of systemic hypertension and ventricular septal defect. For the latter case, we perform forward uncertainty quantification.
\end{abstract}
\maketitle

\section{Introduction}
Forward problems describe mathematical models that map a set of parameters to a corresponding set of observational data. On the other hand, inverse problems aim to determine the model parameters from (noisy) observational data by minimizing a loss function that quantifies the mismatch between the observational data and the model outputs.\\
Deterministic inverse problems seek to identify the single values of the parameters \cite{neto2012introduction}. However, even when the forward map from parameters to outputs exists, the inverse map may not, leading to the ill-posedness of deterministic inverse problems. To address this issue, Tikhonov regularization is often applied to the loss function by adding a term that penalizes the norm of the parameters \cite{tikhonov1977solutions}.\\
Bayesian inverse problems aim to determine the posterior distribution of the parameters, that is their probability distribution given an observational data \cite{stuart2010inverse}. The solution to a Bayesian inverse problem accounts for uncertainty in the parameters estimate that depends on the noise in the observational data.

Neural Networks (NNs) have gained increasing attention for solving inverse problems in recent years \cite{zhu2024label,caforio2024physics,salvador2024digital,salvador2024whole}. They offer a significant reduction in the computational cost of solving inverse problems that would otherwise require numerous iterations of an optimization method \cite{byrd1995limited,chong2013introduction,nelder1965simplex}. A single iteration of the chosen optimization method requires solving the forward problem, which can make inverse problems prohibitively expensive when the forward one is computationally demanding.\\
Tikhonov Networks (TNets) have been shown to effectively solve deterministic inverse problems \cite{nguyen2024tnet}. Their loss function consists of two terms: one incorporating the prior knowledge of the parameters and another enforcing a physics-aware constraint. Data randomization has been theoretically and numerically demonstrated to implicitly regularize the weights and biases of the NN. Tikhonov autoencoders (TAEN) further enhance the use of data randomization, enabling the training of the NN generating new data by the randomization of a single sample \cite{nguyen5081218taen}. InVAErt networks address the ill-posedness nature of deterministic inverse problems by extending the observational data space with latent variables, enabling the identification of different parameters generating the same observational data \cite{tong2024invaert,tong2025invaert}. NETworks Tikhonov (NETT) have been applied to inverse problems in imaging, leveraging a data-driven regularization approach \cite{li2020nett}. Additionally, convolutional Neural Networks (CNN) have been widely explored for image-based inverse problems \cite{daon2018mitigating,jin2017deep,lucas2018using}.\\
In the context of Bayesian inverse problems, Uncertainty Quantification Variational AutoEncoders (UQ-VAE) have been developed to estimate the mean and covariance of the posterior distribution of parameters. In Ref.~\citen{goh2021solving}, the authors established the theoretical foundations of UQ-VAEs, deriving a loss function based on a family of Jensen-Shannon divergences \cite{nielsen2010family}. This approach comes from generative modeling that employs Jensen-Shannon divergences to compute the posterior parameter distribution and generate new data samples \cite{deasy2020constraining,sutter2020multimodal}.  In Ref.~\citen{tonini2025enhanced}, a novel loss function was proposed to improve the estimation accuracy of UQ-VAEs, while enhancing their generalization capabilities. However, this approach required an increased computational cost for evaluating the training loss function, leading to prolonged training times.

In this work, we build upon the results of Ref.~\citen{tonini2025enhanced} proposing a novel loss function by removing a sample mean term from the original loss function, among other modifications. This approach improves accuracy and reduces the training time compared to previous literature results. The time reduction advantage becomes more pronounced as the dimensionality of the parameter space increases. We obtain a theoretical result analogous to the original one for affine forward problems and extend it to general forward problems under additional assumptions, thereby strengthening the theoretical background of UQ-VAEs. We validate our approach on a Poisson problem and compare it in terms of training time and accuracy with the previous one on a non linear problem. Furthermore, we test our approach on two test cases involving 0D cardiocirculatory models \cite{regazzoni2022cardiac,albanese2016integrated,dede2021modeling,tonini2025two,formaggio2025multidomain}, specifically modeling systemic hypertension \cite{world2023global} and ventricular septal defect \cite{spicer2014ventricular,dakkak2024ventricular}. For the latter case, we perform forward uncertainty quantification to further assess the performances of the new approach.

The paper is structured as follows: in \cref{sec:buildingLossFunction}, we introduce the new loss function and UQ-VAEs; in \cref{sec:results}, we present the numerical tests on the Poisson problem and the 0D cardiocirculatory model; in \cref{sec:concl}, we draw our conclusions.

\section{Bayesian inverse problems} \label{sec:buildingLossFunction}
We present the general framework of Bayesian inverse problems and UQ-VAEs. Additionally, we propose a novel loss function for training the UQ-VAEs addressing previous limitations \cite{tonini2025enhanced}. The new loss function reduces the computational cost of training the UQ-VAE while extending the theoretical result of the previous approach.

Bayesian inverse problems are used to estimate the parameters of a mathematical model along with their uncertainties, given noisy observational data. Let $U$ and $Y$ be random variables representing the parameters and the noisy observational data, respectively. We assume the existence of a map $\mathcal F$ from the parameters to the observational data in absence of noise. Let $E$ denote the noise random variable. The mathematical model under consideration is:
\begin{align} \label{eq:model}
Y=\mathcal{F}(U)+E.
\end{align} 
We focus on the case where $U$ and $Y$ are finite dimensional real random variables. Specifically, each parameter is a vector $\mathbf u \in \R^\mathrm{D}$ and each observational data is a vector $\mathbf y \in \R^\mathrm{O}$, with $D,O \in \N$. Thus, the map is defined as $\mathcal F:\R^\mathrm{D} \to \R^\mathrm{O}$. The solution of a Bayesian inverse problem is the posterior distribution of the parameters $p_\mathrm{U|Y}(\mathbf u| \mathbf y)$, which is the distribution of the parameters conditioned to the given noisy observational data.\\
Using the Bayes' theorem \cite{bayes1958essay}, it is possible to retrieve an implicit formula for the posterior distribution that, in most of the cases, is not possible to compute explicitly as it requires the evaluation of an integral. Generally, one looks for the mean and covariance, namely the first two moments, of the posterior distribution to characterize the expected parameter and its uncertainty subjected to the noisy observational data. A standard approach to solve Bayesian inverse problems is by computing the Maximum A Posteriori estimate (MAP) parameter $\mathbf u_\mathrm{MAP} \in \R^\mathrm{D}$, that is the most likely parameter to observe given a noisy observational data $\mathbf y \in \R^\mathrm{O}$, and the Laplace approximation of the covariance matrix $\Gamma_\mathrm{Lap} \in \R^\mathrm{D \times D}$\cite{evans2000approximating}. By the Bayes' theorem, the posterior distribution can be expressed as:
\begin{align}
p_\mathrm{U|Y}(\mathbf u | \mathbf y) = \frac{p_\mathrm{U,Y}(\mathbf u, \mathbf y)}{p_\mathrm{Y}(\mathbf y)} = \frac{p_\mathrm{Y|U}(\mathbf y | \mathbf u)\,p_\mathrm{U}(\mathbf u)}{p_\mathrm{Y}(\mathbf y)} = \frac{p_\mathrm{Y|U}(\mathbf y | \mathbf u)\,p_\mathrm{U}(\mathbf u)}{\int_{\R^\mathrm{D}} p_\mathrm{Y|U}(\mathbf y|\mathbf u)p_\mathrm{U}(\mathbf u) \mathrm{d}\mathbf u}, \label{eq:postProb}
\end{align}
where $p_\mathrm{U,Y}(\mathbf u, \mathbf y)$, $p_\mathrm{Y}(\mathbf y)$, $p_\mathrm{U}(\mathbf u)$ and $p_\mathrm{Y|U}(\mathbf y | \mathbf u)$ are the probability density functions (pdf) of the joint probability distribution of $U$ and $Y$, the marginal pdfs of $Y$ and $U$, and the conditioned probability of $Y$ with respect to $U$ (likelihood), respectively. We assume that $E$ and $U$ are independent, so the likelihood function is:
\begin{align}
p_\mathrm{Y|U}(\mathbf y| \mathbf u) = p_\mathrm{E}(\mathbf y- \mathcal F(\mathbf u)), \label{eq:noiseProb}
\end{align}
where $p_\mathrm{E}$ is the pdf of the noise random variable $E$. Furthermore, we assume that $U\sim \mathcal N(\boldsymbol \mu_\mathrm{pr}, \Gamma_\mathrm{pr})$ and $E \sim \mathcal N(\boldsymbol \mu_\mathrm{E}, \Gamma_\mathrm{E})$, where $\boldsymbol \mu_\mathrm{pr} \in \R^\mathrm{D}$ and $\boldsymbol \mu_\mathrm{E} \in \R^\mathrm{O}$ are the means of the parameters and noise random variables, respectively, and $\Gamma_\mathrm{pr}\in \R^\mathrm{D \times D}$ and $\Gamma_\mathrm{E} \in \R^\mathrm{O\times O}$ are the corresponding covariance matrices. In this work, we restrict our analysis to symmetric positive definite (SPD) covariance matrices, as well as their inverses, to ensure the well definition of the loss function used in training the variational autoencoders. Using \eqref{eq:noiseProb} and the normality assumptions on $U$ and $E$, the joint pdf $p_\mathrm{U,Y}(\mathbf u, \mathbf y)$ can be rewritten as:
\begin{align}
p_\mathrm{U,Y}(\mathbf u , \mathbf y) = p_\mathrm{E}(\mathbf y- \mathcal F(\mathbf u))\,&p_\mathrm{U}(\mathbf u) = \frac{1}{(2\pi)^\mathrm{\frac{O+D}{2}}\lvert \Gamma_\mathrm{E} \rvert^\mathrm{\frac{1}{2}} \lvert \Gamma_\mathrm{pr} \rvert^\mathrm{\frac{1}{2}}}\\
&\exp \left ( -\frac{1}{2}\left( \norm{\mathbf y-\mathcal F( \mathbf u) - \boldsymbol \mu_\mathrm{E}}_\mathrm{\Gamma_\mathrm{E}^\mathrm{-1}}^\mathrm{2}+\norm{\mathbf u -\boldsymbol \mu_\mathrm{pr}}_\mathrm{\Gamma_\mathrm{pr}^\mathrm{-1}}^\mathrm{2} \right) \right),
\end{align}
where $\lvert \cdot \rvert$ denotes the matrix determinant and $\norm{\cdot}_\mathrm{\Gamma^\mathrm{-1}}$ represents the norm induced by an SPD matrix $\Gamma^\mathrm{-1}$. Maximizing the posterior pdf $p_\mathrm{U|Y}(\mathbf u|\mathbf y)$ with respect to $\mathbf u$ yields the MAP estimate:
\begin{align}
\mathbf u_\mathrm{MAP} &\coloneqq \underset{\mathbf u\in \R^\mathrm{D}}{\mathrm{argmax}}\, p_\mathrm{U|Y}(\mathbf u | \mathbf y) \\
&= \underset{\mathbf u\in \R^\mathrm{D}}{\mathrm{argmin}}\left( \norm{\mathbf y-\mathcal F(\mathbf u) -\boldsymbol \mu_\mathrm{E}}^\mathrm{2}_\mathrm{\Gamma_\mathrm{E}^\mathrm{-1}} +  \norm{\mathbf u -\boldsymbol \mu_\mathrm{pr}}^\mathrm{2}_\mathrm{\Gamma_\mathrm{pr}^\mathrm{-1}}\right).
\label{eq:mapMu}
\end{align}
The first term of the optimization problem derives from the likelihood function and quantifies the mismatch between the observational data and model predictions. The second term encodes prior knowledge about the parameters and serves as a regularization term. Solving this optimization problem typically requires gradient-based methods, which can be computationally expensive, particularly when evaluating the map $\mathcal F$ is costly.\\
To quantify uncertainty in the parameter estimate due to noise, the Laplace approximation $\Gamma_\mathrm{Lap}$ of the covariance matrix can be used \cite{evans2000approximating}:
\begin{align}
\Gamma_\mathrm{Lap} = \left(J_\mathrm{\mathcal F}(\mathbf u_\mathrm{MAP})^\mathrm{T} \Gamma_\mathrm{E}^\mathrm{-1}J_\mathrm{\mathcal F}(\mathbf u_\mathrm{MAP})+\Gamma_\mathrm{pr}^\mathrm{-1} \right )^\mathrm{-1}, \label{eq:Lap}
\end{align}
where $J_\mathrm{\mathcal F}(\mathbf u_\mathrm{MAP})\in \R^\mathrm{O\times D}$ is the Jacobian of $\mathcal F$ evaluated at $\mathbf u_\mathrm{MAP}$. Computing $\Gamma_\mathrm{Lap}$ requires evaluating $J_\mathrm{\mathcal F}(\mathbf u_\mathrm{MAP})$, which is significantly more expensive than a single evaluation of $\mathcal F$.

When evaluating $\mathcal F$ is expensive (such as when $\mathcal F$ represents the solution operator of a complex partial differential equation), computing $\mathbf u_\mathrm{MAP}$ and $\Gamma_\mathrm{Lap}$ becomes infeasible. UQ-VAEs \cite{goh2021solving,tonini2025enhanced}, a class of NNs, offer a way to mitigate this computational burden. Although the most advanced UQ-VAE formulation \cite{tonini2025enhanced} presents strong generalization capabilities, even with small datasets, it has two main drawbacks when $\mathcal F$ is not affine: an inaccurate estimate of the posterior covariance matrix and a high computational cost for training. Even if the latter is compensated by the efficient solution of Bayesian inverse problems, reducing the training cost is a priority. In this work, we propose a modification to the loss function used for training UQ-VAEs \cite{tonini2025enhanced}. The new loss function extends  the previous theoretical result while significantly reducing computational costs. This reduction enables also faster testing of the NN architectures, indirectly improving the accuracy of the posterior mean and covariance matrix estimates.

Employing variational inference, we approximate the posterior distribution $p_\mathrm{U|Y}(\mathbf u|\mathbf y)$ using a Gaussian distribution:
\begin{align} 
q_\mathrm{\phi}(\mathbf u|\mathbf y) = \mathcal N(\boldsymbol \mu_\mathrm{post}(\mathbf y, \phi),\Gamma_\mathrm{post}(\mathbf y, \phi)), \label{eq:vaeApp}
\end{align}
where $\boldsymbol \mu_\mathrm{post}(\mathbf y, \phi)\in \R^\mathrm{D}$,  $\Gamma_\mathrm{post}(\mathbf y, \phi)\in \R^\mathrm{D\times D}$ and $\phi$ denotes a set of hyper-parameters (in \cref{sec:UQ-VAE}, $\phi$ is the set of weights and biases of the UQ-VAE). For conciseness, we omit the explicit dependence of $\boldsymbol \mu_\mathrm{post}$ and $\Gamma_\mathrm{post}$ on $\mathbf y$ and $\phi$.

\subsection{Loss function}
We describe the loss function and the theoretical result presented in Ref.~\citen{tonini2025enhanced}, modifying the former to reduce its computational cost.\\
Given $\alpha \in (0,1)$ and an observational data $\mathbf y\in \R^\mathrm{O}$, the loss function proposed in the previous work is:
\begin{gather}
L_\mathrm{\alpha}(\boldsymbol \mu, C) = (1-\alpha) \left( \norm{\boldsymbol \mu - \boldsymbol \mu_\mathrm{pr}}_\mathrm{\Gamma^\mathrm{-1}}^\mathrm{2} + \mathrm{tr}\left(\Gamma^\mathrm{-1}\Gamma_\mathrm{pr}  \right) \right) +\\
\alpha \left( \norm{\boldsymbol \mu - \boldsymbol \mu_\mathrm{pr}}_\mathrm{\Gamma_\mathrm{pr}^\mathrm{-1}}^\mathrm{2} +\mathrm{tr}\left( \Gamma_\mathrm{pr}^\mathrm{-1} \Gamma  \right)  \right)+ \\
 \alpha \E_\mathrm{\mathcal N(\boldsymbol \mu, \Gamma)}\left[ \norm{\mathbf y - \boldsymbol \mu_\mathrm{E} - \mathcal F(\mathbf u)}_\mathrm{\Gamma_\mathrm{E}^\mathrm{-1}}^\mathrm{2} \right], \label{eq:minProblem}
\end{gather}
where $\E_\mathrm{\mathcal N(\boldsymbol \mu, \Gamma)}\left[\cdot \right]$ is the expected value with respect to the probability distribution $\mathcal N(\boldsymbol \mu, \Gamma)$, $\boldsymbol \mu \in \R^\mathrm{D}$, $\Gamma \in \R^\mathrm{D \times D}$ is SPD and $C$ is its Cholesky factor (a lower triangular matrix with positive diagonal entries). This loss function is an upper bound for a family of Jensen-Shannon divergences (JSD) \cite{nielsen2010family} summed to a Kullback-Leibler divergence, which measure the difference between the posterior distribution $p_\mathrm{U|Y}(\mathbf u| \mathbf y)$ and the distribution $\mathcal N(\boldsymbol \mu, \Gamma)$. The first two lines come from the Kullback-Leibler divergences between the prior and the distribution $\mathcal N(\boldsymbol \mu, \Gamma)$ and the third line represents the expected squared error on the observational data.

When $\mathcal F$ is affine (i.e., $\mathcal F(\mathbf u) = F \mathbf u + \mathbf f$, where $F \in \R^\mathrm{O\times D}, \mathbf f \in \R^\mathrm{O}$), \eqref{eq:mapMu} and \eqref{eq:Lap} can be expressed as:
\begin{align}
\mathbf u_\mathrm{MAP} &= \Gamma_\mathrm{Lap} \left( F^\mathrm{T} \Gamma_\mathrm{E}^\mathrm{-1} \left(  \mathbf y - \mathbf f - \boldsymbol \mu_\mathrm{E} \right) +\Gamma_\mathrm{pr}^\mathrm{-1} \boldsymbol \mu_\mathrm{pr}  \right) ,\label{eq:MAPAffine}\\
\Gamma_\mathrm{Lap} &= \left( F^\mathrm{T} \Gamma_\mathrm{E}^\mathrm{-1}F + \Gamma_\mathrm{pr}^\mathrm{-1} \right)^\mathrm{-1}. \label{eq:LapAffine}
\end{align}
Moreover, the expected value in \eqref{eq:minProblem} is computed exactly:
\begin{gather}
 \E_\mathrm{\mathcal N(\boldsymbol \mu, \Gamma)}\left[ \norm{\mathbf y - \boldsymbol \mu_\mathrm{E} - \mathcal F(\mathbf u)}_\mathrm{\Gamma_\mathrm{E}^\mathrm{-1}}^\mathrm{2} \right] =\\
 \norm{\mathbf y - \boldsymbol \mu_\mathrm{E} - \E_\mathrm{\mathcal N(\boldsymbol \mu, \Gamma)}\left[ \mathcal F(\mathbf u) \right]}_\mathrm{\Gamma_\mathrm{E}^\mathrm{-1}}^\mathrm{2} + \mathrm{tr} \left( \Gamma_\mathrm{E}^\mathrm{-1} \mathrm{Cov}_\mathrm{\mathcal N(\boldsymbol \mu, \Gamma)}\left[ \mathcal F(\mathbf u)  \right] \right) = \\
  \norm{\mathbf y - \boldsymbol \mu_\mathrm{E} - F \boldsymbol \mu -\mathbf f }_\mathrm{\Gamma_\mathrm{E}^\mathrm{-1}}^\mathrm{2} + \mathrm{tr} \left( \Gamma_\mathrm{E}^\mathrm{-1} F\Gamma F^\mathrm{T} \right),\label{eq:expValAff}
\end{gather}
where $\mathrm{Cov}_\mathrm{\mathcal N(\boldsymbol \mu, \Gamma)}\left[\cdot \right]$ is the covariance with respect to the probability distribution $\mathcal N(\boldsymbol \mu, \Gamma)$.\\
Minimizing this loss function yields the following result \cite{tonini2025enhanced}.

\begin{thm} \label{thm:convergence}
Assume that $\mathcal F(\mathbf u) = F \mathbf u + \mathbf f$, where $F \in \R^\mathrm{O\times D}, \mathbf f \in \R^\mathrm{O}$, and $\mathbf y \in \R^\mathrm{O}$. Suppose that the prior and noise models are Gaussian and independent, $\mathcal N(\boldsymbol \mu_\mathrm{pr}, \Gamma_\mathrm{pr})$ and $\mathcal N (\boldsymbol \mu_\mathrm{E}, \Gamma_\mathrm{E})$, respectively. Then, the posterior distribution $p_\mathrm{U|Y}(\mathbf u| \mathbf y)$ is $\mathcal N (\mathbf u_\mathrm{MAP}, \Gamma_\mathrm{Lap})$, where $(\mathbf u_\mathrm{MAP}, \Gamma_\mathrm{Lap})$ are given by \eqref{eq:MAPAffine} and \eqref{eq:LapAffine}.\\
Let $\alpha \in \left( 0,1 \right)$ and $(\hat {\boldsymbol \mu}, \hat{C})$ (with $\hat{\Gamma} = \hat{C}\hat{C}^\mathrm{T}$) be the stationary points of the loss function
\begin{align}
 L_\mathrm{\alpha}(\boldsymbol \mu,C)  =&(1-\alpha) \left( \norm{\boldsymbol \mu - \boldsymbol \mu_\mathrm{pr}}_\mathrm{\Gamma^\mathrm{-1}}^\mathrm{2} + \mathrm{tr}\left(\Gamma^\mathrm{-1}\Gamma_\mathrm{pr}  \right) \right) +\\
&\alpha \left( \norm{\boldsymbol \mu - \boldsymbol \mu_\mathrm{pr}}_\mathrm{\Gamma_\mathrm{pr}^\mathrm{-1}}^\mathrm{2} +\mathrm{tr}\left( \Gamma_\mathrm{pr}^\mathrm{-1} \Gamma  \right)  \right) + \\
 & \alpha \left(  \norm{\mathbf y - \boldsymbol \mu_\mathrm{E} - F \boldsymbol \mu -\mathbf f }_\mathrm{\Gamma_\mathrm{E}^\mathrm{-1}}^\mathrm{2} + \mathrm{tr} \left( \Gamma_\mathrm{E}^\mathrm{-1} F\Gamma F^\mathrm{T} \right) \right). \label{eq:minProblemAffine}
\end{align}
Let
\begin{align}
\boldsymbol \mu_\mathrm{post} &= \frac{1-\alpha}{\alpha} \Gamma_\mathrm{Lap} \hat{\Gamma}^\mathrm{-1} (\hat{\boldsymbol \mu}-\boldsymbol \mu_\mathrm{pr}) + \hat{\boldsymbol \mu}, \label{eq:muMap2}\\
\Gamma_\mathrm{post} &= \hat{\Gamma}A^\mathrm{-1} \hat{\Gamma}, \label{eq:gammaMap2}
\end{align}
with
\begin{align}
A =  \frac{1-\alpha}{\alpha}\left( (\hat{\boldsymbol \mu}- \boldsymbol \mu_\mathrm{pr})(\hat{\boldsymbol \mu}- \boldsymbol \mu_\mathrm{pr})^\mathrm{T} + \Gamma_\mathrm{pr} \right). \label{eq:Amatrix2}
\end{align}
Then $(\boldsymbol \mu_\mathrm{post}, \Gamma_\mathrm{post})=(\mathbf u_\mathrm{MAP},\Gamma_\mathrm{Lap})$.
\end{thm}

\emph{Remark.} In practice, we minimize the loss function \eqref{eq:minProblemAffine} to find $(\hat{\boldsymbol \mu},\hat{\Gamma})$ and we compute $(\boldsymbol \mu_\mathrm{post}, \Gamma_\mathrm{post})$ using \eqref{eq:muMap2} and \eqref{eq:gammaMap2}. \\
\emph{Remark.} The matrix $A$ is SPD since it is a positive scalar multiple of the sum of two symmetric matrices that are positive semidefinite and positive definite. Consequently, $\hat \Gamma A^\mathrm{-1} \hat \Gamma$ is also SPD, as $\Gamma_\mathrm{post}$ and $\Gamma_\mathrm{Lap}$ are.

When the map $\mathcal F$ is not affine, $\E_\mathrm{\mathcal N(\boldsymbol \mu, \Gamma)}\left[ \norm{\mathbf y - \boldsymbol \mu_\mathrm{E} - \mathcal F(\mathbf u)}_\mathrm{\Gamma_\mathrm{E}^\mathrm{-1}}^\mathrm{2} \right]$ cannot be computed in a closed form. To approximate this expectation, the central limit theorem is used, leading to the following loss function:
\begin{gather}
 (1-\alpha) \left( \norm{\boldsymbol \mu - \boldsymbol \mu_\mathrm{pr}}_\mathrm{\Gamma^\mathrm{-1}}^\mathrm{2} + \mathrm{tr}\left(\Gamma^\mathrm{-1}\Gamma_\mathrm{pr}  \right) \right) +\\
 \alpha \left(\norm{\boldsymbol \mu - \boldsymbol \mu_\mathrm{pr}}_\mathrm{\Gamma_\mathrm{pr}^\mathrm{-1}}^\mathrm{2} +\mathrm{tr}\left( \Gamma_\mathrm{pr}^\mathrm{-1} \Gamma  \right)  \right)+ \\
 \alpha \frac{1}{K} \sum_\mathrm{k = 1}^\mathrm{K} \norm{\mathbf y - \boldsymbol \mu_\mathrm{E} - \mathcal F(\mathbf u^\mathrm{k}_\mathrm{draw})}_\mathrm{\Gamma_\mathrm{E}^\mathrm{-1}}^\mathrm{2}, \label{eq:minProblem2}
\end{gather}
where $K\in \N$ and, for $k = 1, \dots, K$, $\mathbf u^\mathrm{k}_\mathrm{draw}= \boldsymbol \mu+C \boldsymbol  \epsilon^\mathrm{k}$, with $\boldsymbol  \epsilon^\mathrm{k} \in \R^\mathrm{D}$ and $\boldsymbol  \epsilon^\mathrm{k} \sim \mathcal N(\mathbf 0,I)$. Therefore, at each iteration of the optimization algorithm, $\mathbf u^\mathrm{k}_\mathrm{draw}$ are drawn from $\mathcal N(\boldsymbol \mu, \Gamma)$.

Approximating the expectation in \eqref{eq:minProblem} is typically expensive, even when Quasi-Monte Carlo methods are used \cite{hung2024review}, as a large number of samples $K$ is required for convergence. Given the Cholesky factorizations of $\Gamma_\mathrm{pr}$ and $\Gamma_\mathrm{E}$ the computational cost of evaluating \eqref{eq:minProblem2} is $O(D^\mathrm{3}+KD^\mathrm{2}+K\mathrm{c}(\mathcal F) + KO^\mathrm{2})$, where $\mathrm{c}(\mathcal F)$ is the cost of an evaluation of $\mathcal F$. $O(D^\mathrm{3})$ is the cost of solving $D$ linear systems with matrix $\Gamma_\mathrm{post}$ or $\Gamma_\mathrm{pr}$. $O(KD^\mathrm{2})$, $O(K\mathrm{c}(\mathcal F))$ and $O(KO^\mathrm{2})$ are the costs of computing $K$ samples $\mathbf u_\mathrm{draw}^\mathrm{k}$, of $K$ evaluations of $\mathcal F$ and of solving $K$ linear systems with matrix $\Gamma_\mathrm{E}$, respectively. The computational costs of the other terms are lower than the previous ones.

We propose a novel loss function that significantly reduces the previous computational cost:
\begin{gather}
L_\mathrm{\theta}(\boldsymbol \mu,C) =\theta^\mathrm{2} \mathrm{tr}\left(\Gamma^\mathrm{-1}\Gamma_\mathrm{pr}  \right) + 
\norm{\boldsymbol \mu - \boldsymbol \mu_\mathrm{pr}}_\mathrm{\Gamma_\mathrm{pr}^\mathrm{-1}}^\mathrm{2} +\theta^\mathrm{2} \mathrm{tr}\left( \Gamma_\mathrm{pr}^\mathrm{-1} \Gamma  \right) +\\
\norm{\mathbf y -\boldsymbol \mu_\mathrm{E}-\mathcal F(\boldsymbol \mu)}_\mathrm{\Gamma_\mathrm{E}^\mathrm{-1}}^\mathrm{2} + \\
\mathrm{tr}\left(\Gamma_\mathrm{E}^\mathrm{-1}\left(\tilde {\mathcal F}(\theta C+\boldsymbol \mu \mathbf 1_\mathrm{D}^\mathrm{T}) - \mathcal F (\boldsymbol \mu)\mathbf 1_\mathrm{D}^\mathrm{T} \right)\left(\tilde {\mathcal F}(\theta C+\boldsymbol \mu \mathbf 1_\mathrm{D}^\mathrm{T}) - \mathcal F(\boldsymbol \mu)\mathbf 1_\mathrm{D}^\mathrm{T} \right)^\mathrm{T}  \right),  \label{eq:minProblem3}
\end{gather}
where $\tilde {\mathcal  F}(A) = \left[\mathcal F(A_\mathrm{:,1})\quad \dots \quad \mathcal F(A_\mathrm{:,D})  \right] \in \R^\mathrm{O \times D}$, for $A \in \R^\mathrm{D \times D}$, is the map $\mathcal F$ applied column wise to $A$. The vector $\mathbf 1_\mathrm{D} \in \R^\mathrm{D}$ consists of ones and the optimal value for $\theta \in \R \setminus \{0\}$ is analyzed in \cref{sec:test1}. \\
The new loss function fixes $\alpha = 1/2$ in the previous one \eqref{eq:minProblem} and removes the term $\norm{\boldsymbol \mu - \boldsymbol \mu_\mathrm{pr}}_\mathrm{\Gamma^\mathrm{-1}}^\mathrm{2}$, since it vanishes for $\boldsymbol \mu = \boldsymbol \mu_\mathrm{pr}$ as $\norm{\boldsymbol \mu - \boldsymbol \mu_\mathrm{pr}}_\mathrm{\Gamma_\mathrm{pr}^\mathrm{-1}}^\mathrm{2}$ does. The coefficient $\theta$ is a (small) perturbation needed to obtain the following theoretical result. The coefficient $\theta^\mathrm{2}$ balances the derivatives of $L_\mathrm{\theta}$ with respect to $C$. The expected value of \eqref{eq:minProblem} is substituted with the second and third lines of \eqref{eq:minProblem3} which are more efficient to compute and mimic the expected value in the affine case \eqref{eq:expValAff}.\\
With this novel loss function, we are able to establish the following new theorem.

\begin{thm} \label{thm:convergenceNew}
Assume that $\mathcal F$ is a $C^\mathrm{1}(\R^\mathrm{D},\R^\mathrm{O})$ function, $\mathbf y \in \R^\mathrm{O}$ and that the prior and noise models are Gaussian and independent, $\mathcal N(\boldsymbol \mu_\mathrm{pr}, \Gamma_\mathrm{pr})$ and $\mathcal N (\boldsymbol \mu_\mathrm{E}, \Gamma_\mathrm{E})$, respectively. Let $(\hat {\boldsymbol \mu}(\theta), \hat{C}(\theta))$ be the stationary points of  $L_\mathrm{\theta}(\boldsymbol \mu,C)$. Suppose that $\hat{\underline{\boldsymbol \mu}} = \lim_{\theta \to 0} \hat{\boldsymbol \mu}(\theta)$ and $\hat{\underline{C}} = \lim_{\theta \to 0} \hat{C}(\theta)$ exist and are bounded in norm, $\hat{\underline{C}}$ is invertible and the minimization problem \eqref{eq:mapMu} has only one stationary point $\mathbf u_\mathrm{MAP}$.\\
Let $\hat{\underline{\Gamma}}=\hat{\underline{C}}\hat{\underline{C}}^\mathrm{T}$ and define:
\begin{align}
\boldsymbol \mu_\mathrm{post} &= \hat{\underline{\boldsymbol \mu}}, \label{eq:muMapNew}\\
\Gamma_\mathrm{post} &= \hat{\underline{\Gamma}} \Gamma_\mathrm{pr}^\mathrm{-1} \hat{\underline{\Gamma}}. \label{eq:gammaMapNew}
\end{align}
Then $(\boldsymbol \mu_\mathrm{post}, \Gamma_\mathrm{post})=(\mathbf u_\mathrm{MAP},\Gamma_\mathrm{Lap})$, where $(\mathbf u_\mathrm{MAP},\Gamma_\mathrm{Lap})$ are given by \eqref{eq:mapMu} and \eqref{eq:Lap}.
\end{thm}
\begin{proof}
We compute the gradients of $L_\mathrm{\theta}(\boldsymbol \mu,C)$:
\begin{align}
 \frac{\partial L_\mathrm{\theta}}{\partial \boldsymbol \mu} =& 2\Gamma_\mathrm{pr}^\mathrm{-1}(\boldsymbol \mu-\boldsymbol \mu_\mathrm{pr})-2J_\mathrm{\mathcal F}^\mathrm{T}(\boldsymbol \mu)\Gamma_\mathrm{E}^\mathrm{-1}(\mathbf y-\boldsymbol \mu_\mathrm{E} - \mathcal F (\boldsymbol \mu))+\\
 &2\sum_\mathrm{k = 1}^\mathrm{D}\left(J_\mathrm{\mathcal{F}}(\theta C_\mathrm{:,k}+\boldsymbol \mu)-J_\mathrm{\mathcal{F}}(\boldsymbol \mu) \right)^\mathrm{T}\Gamma_\mathrm{E}^\mathrm{-1}\left(\mathcal{F}(\theta C_\mathrm{:,k}+\boldsymbol \mu)-\mathcal{F}(\boldsymbol \mu) \right),\\
 \frac{\partial L_\mathrm{\theta}}{\partial C} = & -2 \theta^\mathrm{2}\Gamma^\mathrm{-1}\Gamma_\mathrm{pr}C^\mathrm{-T}+2\theta^\mathrm{2}\Gamma_\mathrm{pr}^\mathrm{-1}C+\\
 &2\theta \left[ J^\mathrm{T}_\mathrm{\mathcal{F}}(\theta C_\mathrm{:,1}+\boldsymbol \mu) \Gamma_\mathrm{E}^\mathrm{-1}\left(\mathcal{F}(\theta C_\mathrm{:,1}+\boldsymbol \mu)-\mathcal{F}(\boldsymbol \mu) \right) \quad \dots \right].
\end{align}
Setting the gradients to vanish and omitting the dependence of $(\hat {\boldsymbol \mu}(\theta), \hat{C}(\theta))$ on $\theta$, we get:
\begin{align}
 & \Gamma_\mathrm{pr}^\mathrm{-1}(\hat{\boldsymbol \mu}-\boldsymbol \mu_\mathrm{pr})-J_\mathrm{\mathcal F}^\mathrm{T}(\hat{\boldsymbol \mu})\Gamma_\mathrm{E}^\mathrm{-1}(\mathbf y-\boldsymbol \mu_\mathrm{E} - \mathcal F (\hat{\boldsymbol \mu}))+\\
&\qquad  \sum_\mathrm{k = 1}^\mathrm{D}\left(J_\mathrm{\mathcal{F}}(\theta \hat{C}_\mathrm{:,k}+\hat{\boldsymbol \mu})-J_\mathrm{\mathcal{F}}(\hat{\boldsymbol \mu}) \right)^\mathrm{T}\Gamma_\mathrm{E}^\mathrm{-1}\left(\mathcal{F}(\theta \hat{C}_\mathrm{:,k}+\hat{\boldsymbol \mu})-\mathcal{F}(\hat{\boldsymbol \mu}) \right) = \mathbf{0}, \label{eq:minMu}\\
&-\hat{\Gamma}^\mathrm{-1}\Gamma_\mathrm{pr}\hat{\Gamma}^\mathrm{-1}+\Gamma_\mathrm{pr}^\mathrm{-1}+\\
&\qquad \frac{1}{\theta} \left[ J^\mathrm{T}_\mathrm{\mathcal{F}}(\theta \hat{C}_\mathrm{:,1}+\hat{\boldsymbol \mu}) \Gamma_\mathrm{E}^\mathrm{-1}\left(\mathcal{F}(\theta \hat{C}_\mathrm{:,1}+\hat{\boldsymbol \mu})-\mathcal{F}(\hat{\boldsymbol \mu}) \right) \quad \dots \right] \hat{C}^\mathrm{-1} = 0. \label{eq:minC}
\end{align}
As $\mathcal F \in C^\mathrm{1}(\R^\mathrm{D},\R^\mathrm{O})$ and $\hat{\underline{C}}$ is bounded in norm, we can expand the terms $\mathcal{F}(\theta \hat{C}_\mathrm{:,k}+\hat{\boldsymbol \mu})$ and $J_\mathrm{\mathcal{F}}(\theta \hat{C}_\mathrm{:,k}+\hat{\boldsymbol \mu})$ for $\theta \to 0$ as:
\begin{gather}
 \mathcal{F}(\theta \hat{C}_\mathrm{:,k}+\hat{\boldsymbol \mu}) = \mathcal{F}(\hat{\boldsymbol \mu})+\theta J_\mathrm{\mathcal{F}}(\hat{\boldsymbol \mu}) \hat{C}_\mathrm{:,k} + o(\theta \norm{\hat{C}_\mathrm{:,k}}_\mathrm{2}) = \mathcal{F}(\hat{\boldsymbol \mu})+\theta J_\mathrm{\mathcal{F}}(\hat{\boldsymbol \mu}) \hat{C}_\mathrm{:,k} + o(\theta),\\
 J_\mathrm{\mathcal{F}}(\theta \hat{C}_\mathrm{:,k}+\hat{\boldsymbol \mu}) = J_\mathrm{\mathcal{F}}(\hat{\boldsymbol \mu}) + o(\norm{\hat{C}_\mathrm{:,k}}_\mathrm{2}) = J_\mathrm{\mathcal{F}}(\hat{\boldsymbol \mu}) + o(1),
\end{gather}
where $o(\cdot)$ indicates the small-o notation. The boundedness of $\hat{\underline{C}}$ is needed to remove $\norm{\hat{C}_\mathrm{:,k}}_\mathrm{2}$ from $o(\cdot)$. By computing the limit for $\theta \to 0$ of equation \eqref{eq:minMu}, using $\mathcal F \in C^\mathrm{1}(\R^\mathrm{D},\R^\mathrm{O})$ and that $\hat{\underline{\boldsymbol \mu}}, \hat{\underline{C}}$ exist and are bounded, we obtain:
\begin{align}
 & \Gamma_\mathrm{pr}^\mathrm{-1}(\hat{\boldsymbol \mu}-\boldsymbol \mu_\mathrm{pr})-J_\mathrm{\mathcal F}^\mathrm{T}(\hat{\boldsymbol \mu})\Gamma_\mathrm{E}^\mathrm{-1}(\mathbf y-\boldsymbol \mu_\mathrm{E} - \mathcal F (\hat{\boldsymbol \mu}))+o(\theta) \to \\
&\Gamma_\mathrm{pr}^\mathrm{-1}\left( \hat{\underline{\boldsymbol \mu}}-\boldsymbol \mu_\mathrm{pr} \right)-J_\mathrm{\mathcal F}^\mathrm{T}\left( \hat{\underline{\boldsymbol \mu}}\right)\Gamma_\mathrm{E}^\mathrm{-1}\left( \mathbf y-\boldsymbol \mu_\mathrm{E} - \mathcal F \left( \hat{\underline{\boldsymbol \mu}}\right)\right) =  \mathbf{0},\label{eq:minMu2}
\end{align}
By the above equation, we observe that $\boldsymbol \mu_\mathrm{post} = \hat{\underline{\boldsymbol \mu}}$ satisfies the same equation of the stationary point $\mathbf u_\mathrm{MAP}$ of equation \eqref{eq:mapMu}, therefore, by the uniqueness of the stationary point of \eqref{eq:mapMu}, $\boldsymbol \mu_\mathrm{post}$ is equal to $\mathbf u_\mathrm{MAP}$.\\
Since $\hat{\underline{C}}$ is invertible, $\hat{\underline{\Gamma}}^\mathrm{-1}$ exists. We compute the limit of equation \eqref{eq:minC} for $\theta \to 0$:
\begin{align}
& -\hat{\Gamma}^\mathrm{-1}\Gamma_\mathrm{pr}\hat{\Gamma}^\mathrm{-1}+\Gamma_\mathrm{pr}^\mathrm{-1}+\frac{1}{\theta} \left(J^\mathrm{T}_\mathrm{\mathcal{F}}(\hat{\boldsymbol \mu}) + o(1)\right) \Gamma_\mathrm{E}^\mathrm{-1} \left( \theta J_\mathrm{\mathcal{F}}(\hat{\boldsymbol \mu}) \hat C + o(\theta) \right)\hat{C}^\mathrm{-1} \to \\
&-\hat{\underline{\Gamma}}^\mathrm{-1}\Gamma_\mathrm{pr} \hat{\underline{\Gamma}}^\mathrm{-1}+\Gamma_\mathrm{pr}^\mathrm{-1}+ J^\mathrm{T}_\mathrm{\mathcal{F}}\left( \hat{\underline{\boldsymbol \mu}} \right) \Gamma_\mathrm{E}^\mathrm{-1}J_\mathrm{\mathcal{F}}\left( \hat{\underline{\boldsymbol \mu}} \right) =0. \label{eq:minC2}
\end{align}
Since $ \hat{\underline{\boldsymbol \mu}}$ is equal to $\mathbf u_\mathrm{MAP}$, we get:
\begin{align}
\Gamma_\mathrm{pr}^\mathrm{-1}+ J^\mathrm{T}_\mathrm{\mathcal{F}}\left( \hat{\underline{\boldsymbol \mu}} \right) \Gamma_\mathrm{E}^\mathrm{-1}J^\mathrm{T}_\mathrm{\mathcal{F}}\left( \hat{\underline{\boldsymbol \mu}} \right) = \Gamma_\mathrm{Lap}^\mathrm{-1}.
\end{align}
Therefore, from equation \eqref{eq:minC2}, $\Gamma_\mathrm{post}$ is equal to $\Gamma_\mathrm{Lap}$.
\end{proof}

\emph{Remark.} The above theorem only holds for $\theta \to 0$. In practice, $\theta$ can be chosen arbitrarily small without incurring additional computational cost, as long as machine precision does not hinder the evaluation of the loss.\\
\emph{Remark.} The $\lim_\mathrm{\theta \to 0} L_\mathrm{\theta}(\boldsymbol \mu,C)$ corresponds to the loss function of minimization problem \eqref{eq:mapMu}, which makes the result $\boldsymbol \mu_\mathrm{post} = \mathbf u_\mathrm{MAP}$ intuitive. However, it must be shown that the minimization of $L_\mathrm{\theta}(\boldsymbol \mu,C)$ with respect to $\boldsymbol \mu$ and the limit as $\theta \to 0$ commute, as we have done. Moreover, although taking $\theta \to 0$ recovers the minimization problem for $\mathbf u_\mathrm{MAP}$, the $\theta$ dependent terms still enable the estimation of the posterior covariance matrix.\\
\emph{Remark.} When $\mathcal F(\mathbf u) = F \mathbf u + \mathbf f$ is affine, the result of \Cref{thm:convergenceNew} holds $\forall \, \theta \in \R \setminus \{0\}$ without any hypotheses on $\hat{\underline{\boldsymbol \mu}}$, $\hat{\underline{C}}$ and the minimization problem \eqref{eq:mapMu}. It is sufficient to substitute in \eqref{eq:minMu} and \eqref{eq:minC} the expression of $\mathcal F$ and $J_\mathrm{\mathcal F}(\mathbf u) = F$.\\
\emph{Remark.} The above theorem does not assume a specific posterior distribution therefore the mean and covariance estimates hold in general. We choose to approximate the posterior distribution using a multivariate Gaussian distribution, as shown in \eqref{eq:vaeApp}.\\
\emph{Remark.} The assumption of a unique stationary point for \eqref{eq:mapMu} is restrictive and generally holds only when the prior dominates the loss function.

The computational cost of evaluating \eqref{eq:minProblem3} is $O(D^\mathrm{3}+D\mathrm{c}(\mathcal F) + DO^\mathrm{2})$, which follows from the same cost analysis as for \eqref{eq:minProblem2}. Thanks to the curse of dimensionality, the number of samples $K$ to accurately estimate an expected value is generally much higher than the sample space dimension $D$. Consequently, evaluating \eqref{eq:minProblem3} is significantly less expensive than evaluating \eqref{eq:minProblem2}.

\subsection{Uncertainty quantification variational autoencoders}\label{sec:UQ-VAE}
\begin{figure}[t!]
	\centering
 	\includegraphics[width=\linewidth, height=4.5cm,keepaspectratio]{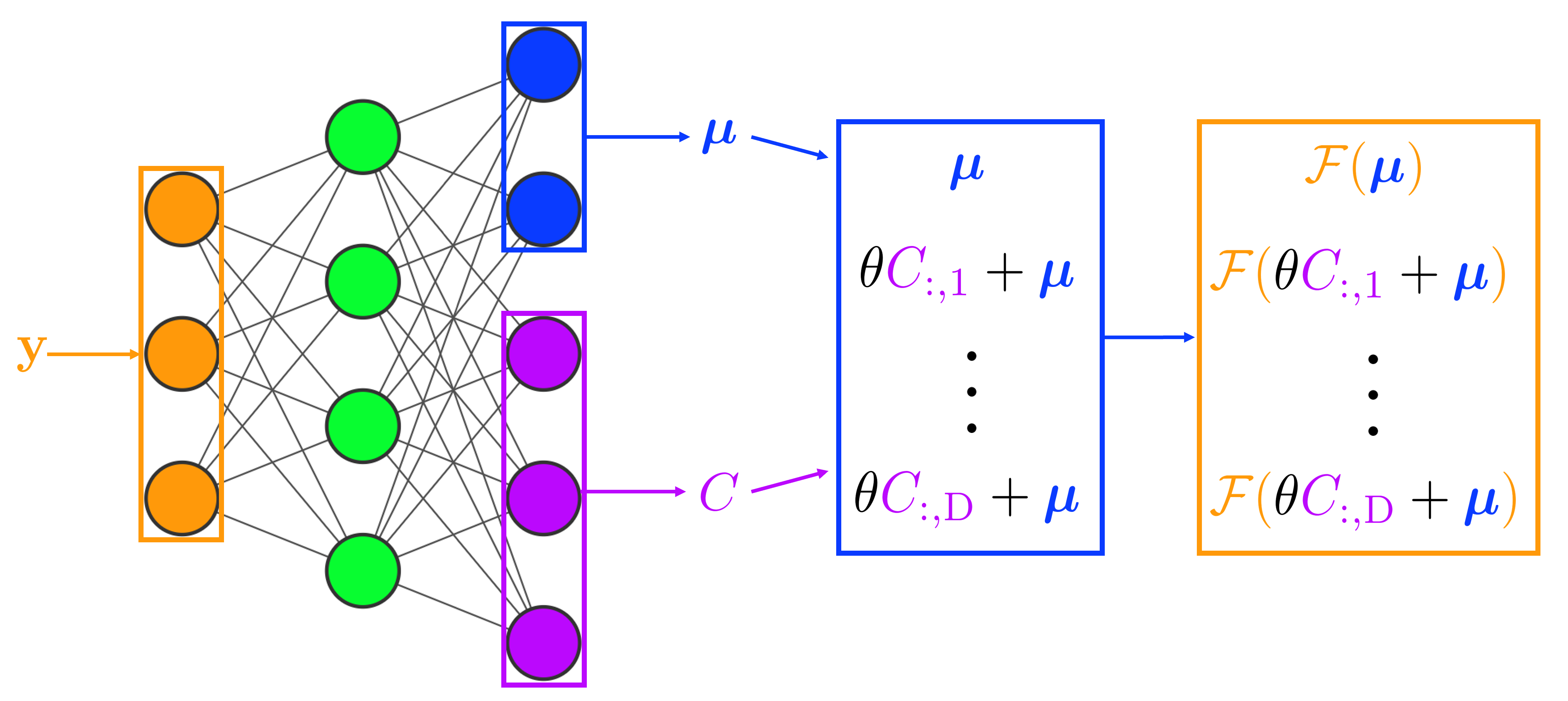}
	\caption{Structure of the UQ-VAE when $\mathcal F$ is evaluated.}
	\label{fig:VAEF}
\end{figure}

An autoencoder consists in the composition of two maps called encoder and decoder. The encoder maps the inputs of the autoencoder to a latent space, while the decoder maps the outputs of the encoder to the original space. Autoencoders typically approximate the identity function, meaning that their inputs and outputs should match. In our method, the encoder is a NN $\varphi: \R^\mathrm{O} \to \R^{D+\frac{D(D+1)}{2}}$, while the decoder corresponds to the map $\mathcal F$ or an approximation of it. The latent space of the UQ-VAE represents the mean $\boldsymbol \mu$ and the Cholesky factor $C$ of the multivariate normal distribution $\mathcal N(\boldsymbol \mu, \Gamma)$ (\Cref{fig:VAEF}).\\
The following Corollary, proven in Ref.~\citen{tonini2025enhanced} and easy to extend to the new loss function \eqref{eq:minProblem3}, demonstrates the convergence of the outputs of a single-layer linear encoder to the posterior mean and covariance in the case of $\mathcal F$ affine.
\begin{cor} \label{cor:convergenceNN}
Consider a single-layer linear encoder. Let $L_\mathrm{\theta}(\phi) = L_\mathrm{\theta}(\boldsymbol \mu(\phi), \Gamma(\phi))$, where $\phi$ is the set of weights and biases of the NN. Under the assumptions of \Cref{thm:convergenceNew} and $\mathcal F$ affine, we define:
\begin{align}
\boldsymbol \mu &= W_\mathrm{\boldsymbol \mu} \mathbf y + \mathbf b_\mathrm{\boldsymbol \mu}, \label{eq:muPost}\\
C &= \mathrm{vec}_\mathrm{L}^\mathrm{-1}(\mathbf l) + \mathrm{diag}(\boldsymbol \sigma), \label{eq:cPost}\\
\boldsymbol \sigma &= \mathrm{exp} \left(  W_\mathrm{\boldsymbol \sigma} \mathbf y + \mathbf b_\mathrm{\boldsymbol \sigma} \right), \label{eq:sigma}\\
\mathbf l &= W_\mathrm{\mathbf l} \mathbf y + \mathbf b_\mathrm{\mathbf l}, \label{eq:l}
\end{align}
where $W_\mathrm{\boldsymbol \mu} \in \R^\mathrm{D \times O}$ and $W_\mathrm{\boldsymbol \sigma} \in \R^\mathrm{D \times O}, W_\mathrm{\mathbf l} \in \R^\mathrm{\frac{D(D-1)}{2} \times O}$ are the encoder weight matrices, $\mathbf b_\mathrm{\boldsymbol \mu} \in \R^\mathrm{D}, \mathbf b_\mathrm{\boldsymbol \sigma} \in \R^\mathrm{D}$ and $\mathbf b_\mathrm{\mathbf l} \in \R^\mathrm{\frac{D(D-1)}{2}}$ are the encoder biases. The functions $\mathrm{exp}$, $\mathrm{vec}_\mathrm{L}$ and $\mathrm{diag}$ denote the elementwise exponential function, the vectorization of a strictly lower triangular matrix and the construction of a diagonal matrix from its vector diagonal, respectively.\\
The stationary points $\hat \phi$ of $L_\mathrm{\theta}(\phi)$ are such that $(\boldsymbol \mu_\mathrm{post}(\hat \phi), \Gamma_\mathrm{post}(\hat \phi)) = (\mathbf u_\mathrm{MAP},\Gamma_\mathrm{Lap})$.
\end{cor}

For deeper encoders with non linear activation functions, $(\boldsymbol \mu, C)$ is computed as in \eqref{eq:muPost}, \eqref{eq:cPost}, \eqref{eq:sigma} and \eqref{eq:l}, using the encoder's outputs in place of $W_\mathrm{\boldsymbol \mu}\mathbf y + \mathbf b_\mathrm{\boldsymbol \mu}$, $W_\mathrm{\boldsymbol \sigma}\mathbf y + \mathbf b_\mathrm{\boldsymbol \sigma}$ and $W_\mathrm{\mathbf l}\mathbf y + \mathbf b_\mathrm{\mathbf l}$.

\begin{figure}[t!]
	\centering
 	\includegraphics[width=\linewidth, height=6cm,keepaspectratio]{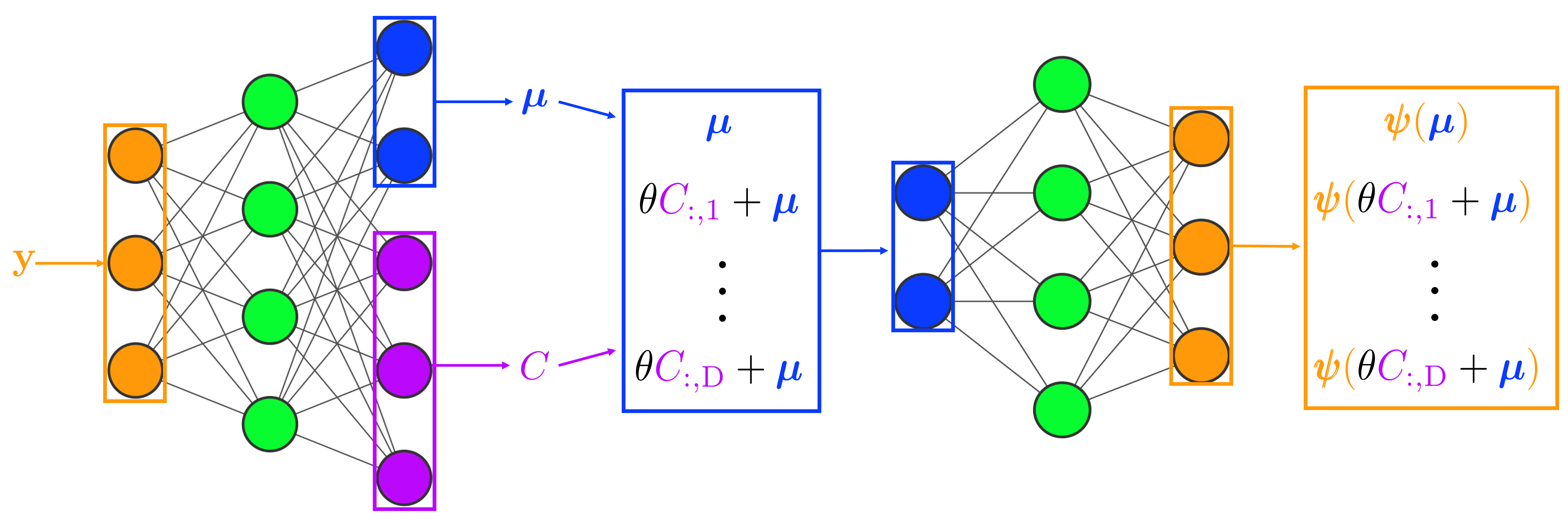}
	\caption{Structure of the UQ-VAE using a NN to approximate $\mathcal F$.}
	\label{fig:VAENN}
\end{figure}

In many applications, evaluating $\mathcal F$ is computationally expensive or $\mathcal F$ can be unknown. In such cases, we approximate the map $\mathcal F$ using a second NN $\psi:\R^\mathrm{D} \to \R^\mathrm{O}$ (\Cref{fig:VAENN}). To account for the approximation error introduced by the NN $\psi$, we use an error random variable $E_\mathrm{dec}\sim \mathcal N(\boldsymbol \mu_\mathrm{dec},\Gamma_\mathrm{dec})$ in the model \eqref{eq:model}:
\begin{align}
\mathcal F(U)=\psi(U)+E_\mathrm{dec} \Rightarrow Y=\psi(U)+E_\mathrm{dec}+E.
\end{align}
We assume that $E_\mathrm{dec}$ and $E$ are independent random variables, so that their sum $E_\mathrm{dec}+E$ is distributed as  $\mathcal N(\boldsymbol \mu_\mathrm{dec}+\boldsymbol \mu_\mathrm{E},\Gamma_\mathrm{dec}+\Gamma_\mathrm{E})$. Defining $\boldsymbol \mu_\mathrm{\tilde E} = \boldsymbol \mu_\mathrm{dec}+\boldsymbol \mu_\mathrm{E}$ and $\Gamma_\mathrm{\tilde E} = \Gamma_\mathrm{dec}+\Gamma_\mathrm{E}$, the corresponding loss function becomes:
\begin{gather}
\theta^\mathrm{2} \mathrm{tr}\left(\Gamma^\mathrm{-1}\Gamma_\mathrm{pr}  \right) +
 \norm{\boldsymbol \mu - \boldsymbol \mu_\mathrm{pr}}_\mathrm{\Gamma_\mathrm{pr}^\mathrm{-1}}^\mathrm{2} +\theta^\mathrm{2} \mathrm{tr}\left( \Gamma_\mathrm{pr}^\mathrm{-1} \Gamma  \right) + \\
  \norm{\mathbf y -\mu_\mathrm{\tilde E}-\psi(\boldsymbol \mu)}_\mathrm{\Gamma_\mathrm{\tilde E}^\mathrm{-1}}^\mathrm{2} +\\ \mathrm{tr}\left(\Gamma_\mathrm{\tilde E}^\mathrm{-1}\left(\tilde {\psi}(\theta C+\boldsymbol \mu \mathbf 1_\mathrm{D}^\mathrm{T}) - \psi (\boldsymbol \mu)\mathbf 1_\mathrm{D}^\mathrm{T} \right)\left(\tilde {\psi}(\theta C+\boldsymbol \mu \mathbf 1_\mathrm{D}^\mathrm{T}) - \psi(\boldsymbol \mu)\mathbf 1_\mathrm{D}^\mathrm{T} \right)^\mathrm{T}  \right),\label{eq:minProblem4}
\end{gather}
where $\tilde \psi(A) = \left[\psi(A_\mathrm{:,1})\quad \dots \quad \psi(A_\mathrm{:,D})  \right] \in \R^\mathrm{O \times D}$, for $A \in \R^\mathrm{D \times D}$.\\
When the UQ-VAE consists of two NNs, its training is performed in two stages: at first the decoder NN is trained and after the encoder one \cite{nguyen5081218taen}. Simultaneous training of both NNs may lead to unpredictable results, as our analysis of the minima of the loss function \eqref{eq:minProblem3} assumes a fixed decoder during optimization.

\emph{Remark.} In practice, we always use fully connected feed forward NNs as encoders and decoders. However, the proposed approach is compatible with other NN architectures. \\
\emph{Remark.} We train the decoder NN using the squared error between its outputs and the noiseless observational data.\\
\emph{Remark.} After training the decoder NN, we estimate $\boldsymbol \mu_\mathrm{dec}$ and $\Gamma_\mathrm{dec}$ by computing the sample mean and covariance of $\mathcal F(U)-\psi(U)$ on a test set of approximately $29,000$ samples.\\
\emph{Remark.} When a dataset $\{\mathbf y^\mathrm{(m)}\}_\mathrm{m=1}^\mathrm{M}$, with $M\in \N$, is available, we train the encoder by minimizing the average of the loss functions \eqref{eq:minProblem3} or \eqref{eq:minProblem4} over the dataset.

\subsubsection{Weights and biases initialization and dataset normalization}\label{sec:WIN}
The initialization of the weights and biases of the encoder plays a crucial role in its training. A good initial guess is obtained by setting $\boldsymbol \mu \sim \boldsymbol \mu_\mathrm{pr}$ and $C \sim C_\mathrm{pr}$ (with $C_\mathrm{pr}$ the Cholesky factor of $\Gamma_\mathrm{pr}$), as these values minimize the term $\theta^\mathrm{2} \mathrm{tr}\left(\Gamma^\mathrm{-1}\Gamma_\mathrm{pr}  \right)+\norm{\boldsymbol \mu - \boldsymbol \mu_\mathrm{pr}}_\mathrm{\Gamma_\mathrm{pr}^\mathrm{-1}}^\mathrm{2} +\theta^\mathrm{2} \mathrm{tr}\left( \Gamma_\mathrm{pr}^\mathrm{-1} \Gamma  \right)$ in the loss function \eqref{eq:minProblem3}. The encoder  weights are initialized according to the Xavier uniform initializer \cite{glorot2010understanding}. Then, all the weights in the last layer are scaled by $10^\mathrm{-4}$. This ensures that, up to a small discrepancy, the output of the encoder is given by the biases of the last layer. All encoder biases are set to zero, except for those of the last layer which are fixed to $\boldsymbol \mu_\mathrm{pr}$ and $\mathbf c_\mathrm{pr}$, where $\mathbf c_\mathrm{pr}$  is given by applying the inverse of the right hand side transformation of \eqref{eq:cPost} to $C_\mathrm{pr} \in \R^\mathrm{D \times D}$.

Since the loss function \eqref{eq:minProblem3} depends on the mean and covariance of both parameters and noise, normalizing the dataset (to enhance the UQ-VAE generalization capability) also affects these quantities. \\
We consider a dataset $\{(\mathbf u^\mathrm{(m)}, \mathbf y^\mathrm{(m)})\}_\mathrm{m=1}^\mathrm{M}$, where $M\in \N$, $\mathbf y^\mathrm{(m)} = \mathcal F(\mathbf u^\mathrm{(m)})+\epsilon^\mathrm{(m)}$, with $\mathbf u^\mathrm{(m)}$ and $\epsilon^\mathrm{(m)}$ are sampled from $\mathcal N(\boldsymbol \mu_\mathrm{pr}, \Gamma_\mathrm{pr})$ and $\mathcal N(\boldsymbol \mu_\mathrm{E}, \Gamma_\mathrm{E})$, respectively. To normalize the data within a given range, e.g., $[0,1]$, we introduce vectors $\mathbf a, \mathbf b \in \R^\mathrm{D}$ and $\mathbf c, \mathbf d\in \R^\mathrm{O}$ such that:
\begin{align}
\bar {\mathbf u}^\mathrm{(m)} = \mathbf u^\mathrm{(m)} \odot \mathbf a+ \mathbf b \in [0,1]^\mathrm{D},\\
\bar {\mathbf y}^\mathrm{(m)} = \mathbf y^\mathrm{(m)} \odot \mathbf c+ \mathbf d \in [0,1]^\mathrm{O}.
\end{align}
As shown in Section $2.4.1$ of Ref.~\citen{tonini2025enhanced}, the normalized random variables $\bar Y$, $\bar U$ and $\bar E$ satisfy a model analogous to \eqref{eq:model} in the normalized scale, with $\bar U \sim N(\bar{\boldsymbol \mu}_\mathrm{pr}, \bar{\Gamma}_\mathrm{pr})$ and $\bar E \sim N(\bar{\boldsymbol \mu}_\mathrm{E}, \bar{\Gamma}_\mathrm{E})$, where:
\begin{align}
\bar{\boldsymbol \mu}_\mathrm{pr} = \boldsymbol \mu_\mathrm{pr} \odot \mathbf a +\mathbf b,\\
\bar{\Gamma}_\mathrm{pr} = \Gamma_\mathrm{pr} \odot (\mathbf a \mathbf a^\mathrm{T}),\\
\bar{\boldsymbol \mu}_\mathrm{E} = \boldsymbol \mu_\mathrm{E} \odot \mathbf c,\\
\bar{\Gamma}_\mathrm{E} = \Gamma_\mathrm{E} \odot (\mathbf c \mathbf c^\mathrm{T}).
\end{align}
The same normalization is applied to observational data both with and without noise. Once $(\bar{\boldsymbol \mu}_\mathrm{post},\bar{\Gamma}_\mathrm{post})$ in the normalized scale is computed, we revert to the original scale using the transformation:
\begin{gather}
\boldsymbol \mu_\mathrm{post} = (\bar{\boldsymbol \mu}_\mathrm{post}-\mathbf b) \oslash \mathbf a,\\
\Gamma_\mathrm{post} = \bar{\Gamma}_\mathrm{post} \odot ((\mathbf 1_\mathrm{D} \oslash \mathbf a)(\mathbf 1_\mathrm{D} \oslash \mathbf a)^\mathrm{T}).
\end{gather}

\section{Numerical results}\label{sec:results}
We carry out three benchmark numerical tests to show the capabilities of UQ-VAEs: in \cref{sec:test1}, we solve Bayesian inverse problems for the Poisson equation and study the impact of varying the number of neurons and layers of both the decoder and the encoder, as well as modifying $\theta$ and the robustness with respect to the noise level $\eta$, the number of observations $O$, the prior distribution $\mathcal N (\boldsymbol \mu_\mathrm{pr},\Gamma_\mathrm{pr})$; in \cref{sec:nonLinear}, we face Bayesian inverse problems involving an exponential map $\mathcal F$, investigate the effect of varying the dimensions of the parameter and data space, $D$ and $O$, and compare the performance, both in terms of computational time and accuracy, of UQ-VAEs trained using the loss function $L_\mathrm{\alpha}$ \eqref{eq:minProblem} or $L_\mathrm{\theta}$ \eqref{eq:minProblem3}; in \cref{sec:cardioTheory}, we address Bayesian inverse problems for a 0D cardiocirculatory model in the context of systemic hypertension (\cref{sec:test2}) and ventricular septal defects (\cref{sec:test3}). The activation functions of the NNs are \emph{ReLU} functions and the data are normalized to the range $[0,1]$. Additionally, we always use the \emph{Adam} optimizer \cite{kingma2014adam}. We tuned the number of neurons and layers of each NN to minimize the validation loss functions. The generated datasets consist of $M$ samples, which are split in $\lfloor 9M/10 \rfloor$ training samples and $M-\lfloor 9M/10 \rfloor$ validation samples, where $\lfloor \cdot \rfloor$ denotes the lower integer part.

In the following sections, the true posterior mean and covariance matrix are approximated using Sobol' sequences \cite{sobol1967distribution}, a Quasi-Monte Carlo method that accelerates the convergence of the sample mean to the expected value. Using \eqref{eq:postProb} and \eqref{eq:noiseProb}, we obtain:
\begin{gather}
\E_\mathrm{p_\mathrm{U|Y}(\mathbf u|\mathbf y)}[\mathbf u] = \frac{\E_\mathrm{p_\mathrm{U}(\mathbf u)}[\mathbf u \, p_\mathrm{E}(\mathbf y- \mathcal F(\mathbf u))]}{\E_\mathrm{p_\mathrm{U}(\mathbf u)}[p_\mathrm{E}(\mathbf y- \mathcal F(\mathbf u))]}, \label{eq:trueMean}\\
\mathrm {Cov}_\mathrm{p_\mathrm{U|Y}[\mathbf u|\mathbf y]}(\mathbf u) = \frac{\E_\mathrm{p_\mathrm{U}(\mathbf u)}[\mathbf u \mathbf u^\mathrm{T}\, p_\mathrm{E}(\mathbf y- \mathcal F(\mathbf u))]}{\E_\mathrm{p_\mathrm{U}(\mathbf u)}[p_\mathrm{E}(\mathbf y- \mathcal F(\mathbf u))]}-\\
\frac{\E_\mathrm{p_\mathrm{U}(\mathbf u)}[\mathbf u \, p_\mathrm{E}(\mathbf y- \mathcal F(\mathbf u))]\E_\mathrm{p_\mathrm{U}(\mathbf u)}[\mathbf u \, p_\mathrm{E}(\mathbf y- \mathcal F(\mathbf u))]^\mathrm{T}}{\E_\mathrm{p_\mathrm{U}(\mathbf u)}[p_\mathrm{E}(\mathbf y- \mathcal F(\mathbf u))]^\mathrm{2}}. \label{eq:trueCov}
\end{gather}
This procedure is computationally expensive and can only be performed for simple problems.

\subsection{Poisson equation} \label{sec:test1}
We consider the Poisson equation on the unit square $\Omega=(0,1)^\mathrm{2}$:
\begin{align} \label{eq:Poisson}
\begin{cases}
-\nabla(e^\mathrm{u(\mathbf x)} \nabla y(\mathbf x)) = \frac{17}{4}\pi^\mathrm{2}\sin(2\pi x_\mathrm{1}) \sin(2\pi x_\mathrm{2}) \quad& \mathbf x \in \Omega,\\
y(\mathbf x) = 0 & \mathbf x \in \partial \Omega,
\end{cases}
\end{align}
where $u(\mathbf x)$ is the (log) diffusion parameter. We approximate the solution $y(\mathbf x)$ of the Poisson problem by the linear finite element method on a triangular mesh of $D = 17^\mathrm{2}=289$ degrees of freedom.\\
We define a Gaussian random field on the whole domain $\Omega$, with mean and covariance functions $\mu_\mathrm{pr} (\mathbf x)$ and $\mathcal C_\mathrm{pr}(\mathbf x,\mathbf z)$, for $\mathbf x, \mathbf z \in (0,1)$. The parameter $U$ is given by the evaluation of this Gaussian random field at the degrees of freedom of the mesh, obtaining $U \sim \mathcal N(\boldsymbol \mu_\mathrm{pr},\Gamma_\mathrm{pr})$. We define $\mu_\mathrm{pr}$ to be zero in $\Omega$ and the covariance function $\mathcal C_\mathrm{pr}$ to ensure the well-posedness of infinite dimensional Bayesian inverse problems \cite{stuart2010inverse,bui2013computational}. We define the differential operator $\mathcal A$:
\begin{align}
\mathcal A u = -\gamma \Delta u+\delta u \quad \text{in } \Omega,
\end{align}
where $\gamma, \delta >0$. The operator $\mathcal A^\mathrm{-2}(\mathbf x, \mathbf z) = G_\mathrm{2}(\mathbf x, \mathbf z)$, for $(\mathbf x, \mathbf z) \in \Omega$, is defined as the solution operator of:
\begin{align}
\begin{cases}
\mathcal A G_\mathrm{1}(\mathbf x, \mathbf z) = \delta_\mathrm{\mathbf x} \qquad &\mathbf z \in \Omega,\\
\nabla G_\mathrm{1}(\mathbf x, \mathbf z) \cdot \mathbf n(\mathbf z) + \beta G_\mathrm{1}(\mathbf x, \mathbf z) = 0  &\mathbf z \in \partial \Omega,\\
\mathcal A G_\mathrm{2}(\mathbf x, \mathbf z) = G_\mathrm{1}(\mathbf x, \mathbf z) & \mathbf z \in \Omega,\\
\nabla G_\mathrm{2}(\mathbf x, \mathbf z) \cdot \mathbf n(\mathbf z) + \beta G_\mathrm{2}(\mathbf x, \mathbf z) = 0 & \mathbf z \in \partial \Omega,
\end{cases} 
\end{align}
where, for $\mathbf x\in \Omega$, $\delta_\mathrm{\mathbf x}$ is the Dirac delta centered in $\mathbf x$, $\mathbf n$ is the unit exterior normal to $\partial \Omega$ and $\beta$ is chosen to reduce the boundary artifacts \cite{daon2018mitigating}. We set $\gamma = 0.1$ and $\delta = 0.5$. The covariance function is given by $\mathcal C_\mathrm{pr}(\mathbf x, \mathbf z) = \mathcal A^\mathrm{-\frac{3}{2}}(\mathbf x,\mathbf z)$. We compute $\mathcal C_\mathrm{pr}$ using the same mesh and finite element method as for solving the Poisson problem, via the Python library for inverse problems \emph{hIPPYlib} \cite{VillaPetraGhattas21,VillaPetraGhattas18,VillaPetraGhattas16}. Finally, we define $(\boldsymbol \mu_\mathrm{pr})_\mathrm{i} = \mu_\mathrm{pr}(\mathbf x_\mathrm{i})$ and $(\Gamma_\mathrm{pr})_\mathrm{i,j} = \mathcal C_\mathrm{pr}(\mathbf x_\mathrm{i},\mathbf x_\mathrm{j})$ for $i,j = 1,\cdots,D$, where $\mathbf x_\mathrm{i}$ are the coordinates of the degrees of freedom of the mesh.\\
The observational data consist of noisy evaluations of the solution $y(\mathbf x)$ of \eqref{eq:Poisson} at $O = 20$ random points in the domain $\Omega$.

We generate a dataset $\{(\mathbf u^\mathrm{(m)}, \mathbf y^\mathrm{(m)}) \}_\mathrm{m=1}^\mathrm{M}$, where $M\in \N$ and $\mathbf u^\mathrm{(m)}$ is sampled from $\mathcal N (\boldsymbol \mu_\mathrm{pr},\Gamma_\mathrm{pr})$. Let $\tilde{\mathbf y}^\mathrm{(m)}$ be the noiseless observations of the solution of \eqref{eq:Poisson} with parameter $\mathbf u^\mathrm{(m)}$. Then, $\mathbf y^\mathrm{(m)} = \tilde{\mathbf y}^\mathrm{(m)}+\boldsymbol \epsilon^\mathrm{(m)}$, for $m = 1,\cdots,M$, where $\boldsymbol \epsilon^\mathrm{(m)}$ is sampled from $\mathcal N(\boldsymbol \mu_\mathrm{E},\Gamma_\mathrm{E})$ with $\boldsymbol \mu_\mathrm{E} = \mathbf 0$ and $\Gamma_\mathrm{E}$ is a diagonal matrix with $(\Gamma_\mathrm{E})_\mathrm{i,i} = (\eta \, \underset{j,m}{\max}\, \abs{\tilde y_\mathrm{j}^\mathrm{(m)}})^\mathrm{2}$, for $i = 1,\cdots, O$ and $\eta \in \R$. We fix $\eta = 0.01$.

To normalize the data to $[0,1]$, we apply a uniform normalization to the parameters and observational data across the entire mesh. We define:
\begin{gather*}
a = (\underset{i,m}{\max} \,u^\mathrm{(m)}_\mathrm{i}-\underset{i,m}{\min} \,u^\mathrm{(m)}_\mathrm{i})^\mathrm{-1},\\
b = -a\, \underset{i,m}{\min} \,u^\mathrm{(m)}_\mathrm{i},\\
c = (\underset{i,m}{\max} \,y^\mathrm{(m)}_\mathrm{i}-\underset{i,m}{\min} \,y^\mathrm{(m)}_\mathrm{i})^\mathrm{-1},\\
d = -c\, \underset{i,m}{\min} \,y^\mathrm{(m)}_\mathrm{i}
\end{gather*}
and set $\bar u(\mathbf x) = u(\mathbf x) a + b$ and $\bar y(\mathbf x) = y(\mathbf x) c + d$ in the continuous framework. Referring to \cref{sec:WIN}, in the discrete setting, $\mathbf a = a\mathbf 1_\mathrm{D}$, $\mathbf b = b\mathbf 1_\mathrm{D}$, $\mathbf c = c\mathbf 1_\mathrm{O}$ and $\mathbf d = d\mathbf 1_\mathrm{O}$. We can derive the Poisson problem in the normalized scale:
\begin{align}
\begin{cases}
-\nabla(e^\mathrm{\frac{\bar u(\mathbf x)}{a}} \nabla \bar y(\mathbf x)) = c e^\mathrm{\frac{b}{a}}\frac{17}{4}\pi^\mathrm{2}\sin(2\pi x_\mathrm{1}) \sin(2\pi x_\mathrm{2}) \quad& \mathbf x \in \Omega,\\
\bar y(\mathbf x) = d & \mathbf x \in \partial \Omega.
\end{cases}
\end{align}

We perform hyperparameter tuning to select the optimal number of layers and neurons per layer for both the decoder and the encoder neural networks. We first determine the best architecture for the decoder network and then use this decoder to train the encoder with parameter $\theta = 10^\mathrm{-4}$. According to the hyperparameter tuning results (\Cref{table:ht}), we select a decoder with $6$ hidden layers and $100$ neurons per layer, trained on a dataset of $M = 4096$ samples. The chosen encoder has $2$ hidden layers and $1000$ neurons per layer and is trained on a subset of $100$ samples of the decoder dataset. Although the dimension of the encoder dataset is small with respect to the decoder one, we show in the following that it still generalizes to unseen data.
\begin{table}[t!]
\scriptsize
\centering
	\begin{tabular}{|c|c|c|c|c|c|c|}
	\hline
	\multicolumn{7}{|c|}{\textbf{Decoder}}\\
          \hline
	Neurons\;\textbackslash \;Layers & $1$ & $2$ & $3$ & $4$ & $5$ & $6$ \\
          \hline
	$50$ & $5.8\cdot 10^\mathrm{-4}$  & $2.8\cdot 10^\mathrm{-4}$  & $2.2\cdot 10^\mathrm{-4}$ & $2.0\cdot 10^\mathrm{-4}$ & $2.3\cdot 10^\mathrm{-4}$ & $2.1\cdot 10^\mathrm{-4}$ \\
	$100$ & $7.4\cdot 10^\mathrm{-4}$  & $3.1\cdot 10^\mathrm{-4}$  & $2.1\cdot 10^\mathrm{-4}$ & $1.8\cdot 10^\mathrm{-4}$ & $1.5\cdot 10^\mathrm{-4}$ & $\mathbf{1.2\cdot 10^\mathrm{-4}}$\\
	$150$ & $6.4\cdot 10^\mathrm{-4}$ & $2.8\cdot 10^\mathrm{-4}$  & $2.3\cdot 10^\mathrm{-4}$ & $1.6\cdot 10^\mathrm{-4}$ & $1.3\cdot 10^\mathrm{-4}$ & $1.4\cdot 10^\mathrm{-4}$\\
	$200$ & $6.4\cdot 10^\mathrm{-4}$  & $3.0\cdot 10^\mathrm{-4}$  & $2.4\cdot 10^\mathrm{-4}$ & $1.8\cdot 10^\mathrm{-4}$ & $2.1\cdot 10^\mathrm{-4}$ & $1.5\cdot 10^\mathrm{-4}$\\
	$250$ & $5.9\cdot 10^\mathrm{-4}$  & $3.4\cdot 10^\mathrm{-4}$  & $2.3\cdot 10^\mathrm{-4}$ & $2.5\cdot 10^\mathrm{-4}$ & $2.1\cdot 10^\mathrm{-4}$ & $1.8\cdot 10^\mathrm{-4}$\\
	\hline
          \end{tabular}\\
	\vspace{1em}
	\begin{tabular}{|c|c|c|c|c|}
	\hline
	\multicolumn{5}{|c|}{\textbf{Encoder}}\\
          \hline
	Neurons\; \textbackslash \; Layers & $1$ & $2$ & $3$ & $4$\\
          \hline
	$500$ & $26.8$  & $24.0$ & $25.2$ & $25.7$  \\
	$750$ & $25.6$ & $24.2$  & $24.6$ & $25.5$ \\
	$1000$ & $24.5$ & $\mathbf{23.5}$ & $24.9$ & $25.6$ \\
	$1250$ & $24.5$ & $24.0$ & $25.3$ & $26.5$ \\
	\hline 
          \end{tabular}
\caption{Final validation loss functions for the decoder and the encoder for the Poisson problem.}
\label{table:ht}
\end{table}

We generate $2^\mathrm{19}$ samples from the distribution $\mathcal N(\boldsymbol \mu_\mathrm{pr},\Gamma_\mathrm{pr})$ to estimate \eqref{eq:trueMean} and \eqref{eq:trueCov}.

We analyze the UQ-VAE performances varying $\theta$ in the loss function \eqref{eq:minProblem4}. For the posterior means $\E_\mathrm{p_\mathrm{U|Y}(\mathbf u|\mathbf y^\mathrm{(m)})}[\mathbf u]$ and variances $\mathrm {Var}_\mathrm{p_\mathrm{U|Y}(\mathbf u|\mathbf y^\mathrm{(m)})}[\mathbf u]$ (the diagonal of \eqref{eq:trueCov}), with $m \in \{1,\dots,M\}$, we compute the relative error with respect to the infinity norms at each degree of freedom:
\begin{gather}
\mathbf {e}_\mathrm{\E_\mathrm{p_\mathrm{U|Y}(\mathbf u|\mathbf y^\mathrm{(m)})}[\mathbf u]} = \frac{\abs{\boldsymbol \mu_\mathrm{post}^\mathrm{(m)}-\E_\mathrm{p_\mathrm{U|Y}(\mathbf u|\mathbf y^\mathrm{(m)})}[\mathbf u]}}{\norm{\E_\mathrm{p_\mathrm{U|Y}(\mathbf u|\mathbf y^\mathrm{(m)})}[\mathbf u]}_\mathrm{\infty}}, \label{eq:relErrVec}\\
\mathbf {e}_\mathrm{\mathrm {Var}_\mathrm{p_\mathrm{U|Y}(\mathbf u|\mathbf y^\mathrm{(m)})}[\mathbf u]} = \frac{\abs{\mathrm{diag}\left(\Gamma_\mathrm{post}^\mathrm{(m)}\right)-\mathrm {Var}_\mathrm{p_\mathrm{U|Y}(\mathbf u|\mathbf y^\mathrm{(m)})}(\mathbf u)}}{\norm{\mathrm {Var}_\mathrm{p_\mathrm{U|Y}(\mathbf u|\mathbf y^\mathrm{(m)})}[\mathbf u]}_\mathrm{\infty}} \label{eq:relErrMat},
\end{gather}
where $(\boldsymbol \mu_\mathrm{post}^\mathrm{(m)},\Gamma_\mathrm{post}^\mathrm{(m)})$ are the UQ-VAE mean and covariance estimates for the $m$-th sample, $\abs{\cdot}$ is applied element-wise and $\norm{\cdot}_\mathrm{\infty}$ is the vector infinite norm. The maximum over the mesh of these errors is the relative error on the $m$-th sample in the infinity norm.
\begin{figure}[t!]
	\centering
 	\includegraphics[width=\linewidth, height=12cm,keepaspectratio]{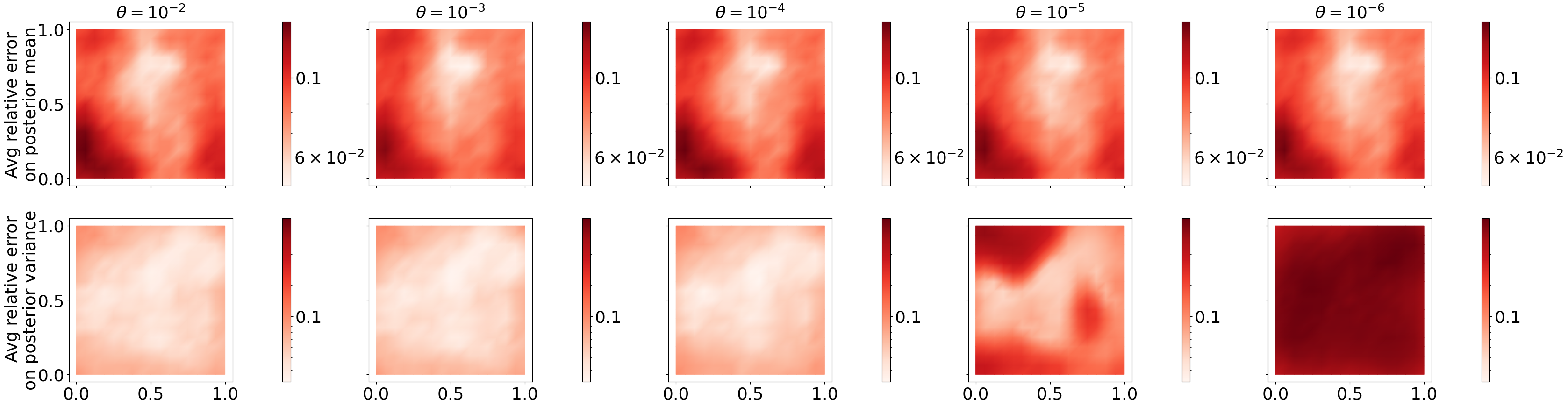}
	\caption{Average relative errors (in logarithmic scale) of the UQ-VAE mean $\boldsymbol \mu_\mathrm{post}^\mathrm{(m)}$ and variance $\Gamma_\mathrm{post}^\mathrm{(m)}$ on the posterior mean $\E_\mathrm{p_\mathrm{U|Y}(\mathbf u|\mathbf y^\mathrm{(m)})}[\mathbf u]$ and on the posterior variance $\mathrm {Var}_\mathrm{p_\mathrm{U|Y}(\mathbf u|\mathbf y^\mathrm{(m)})}[\mathbf u]$ in dependence of $\theta$ for a test set of $100$ samples for the Poisson problem.}
	\label{fig:test1theta}
\end{figure}%

\begin{figure}[t!]
	\centering
 	\includegraphics[width=\linewidth, height=6cm,keepaspectratio]{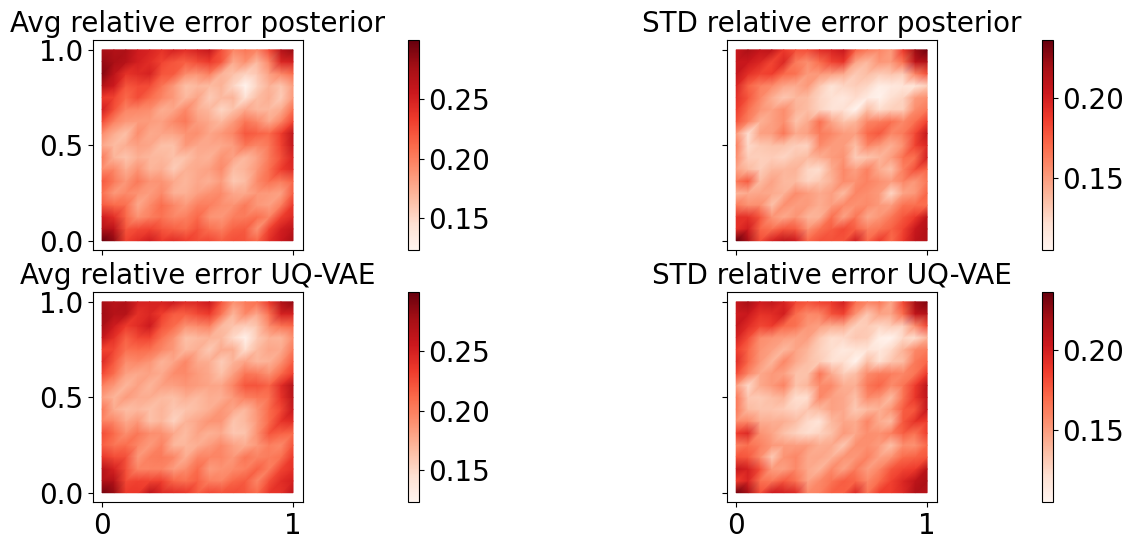}
	\caption{Averages and standard deviations of the relative errors of the posterior mean $\E_\mathrm{p_\mathrm{U|Y}(\mathbf u|\mathbf y^\mathrm{(m)})}[\mathbf u]$ and $\boldsymbol \mu_\mathrm{post}^\mathrm{(m)}$ on the true diffusion coefficients $\mathbf u^\mathrm{(m)}$ for a test set of $100$ samples for the Poisson problem.}
	\label{fig:test1TrueMean}
\end{figure}

The average relative errors on the posterior mean for a test set of $100$ samples share similar trends and magnitudes for all values of $\theta$ (\Cref{fig:test1theta}). The maximum average relative error on the posterior variance is more than $50 \%$ for $\theta=10^\mathrm{-5} \lor 10^\mathrm{-6}$, due to the small gradient of $L_\mathrm{\theta}$ \eqref{eq:minProblem4} with respect to $C$ which scales with $\theta^\mathrm{2}$. For the other values of $\theta$ the errors on the variance are similar. We choose to use in all the tests that follow (also in \cref{sec:nonLinear} and \cref{sec:cardioTheory}) $\theta = 10^\mathrm{-4}$, because \Cref{thm:convergenceNew} holds for $\theta \to 0$.\\
For $\theta = 10^\mathrm{-4}$, the average relative error on the posterior mean and the posterior variance are less than $11\%$. Even if there are differences in $\E_\mathrm{p_\mathrm{U|Y}(\mathbf u|\mathbf y^\mathrm{(m)})}[\mathbf u]$ and $\boldsymbol \mu_\mathrm{post}^\mathrm{(m)}$, they estimate the true diffusion coefficients $\mathbf u^\mathrm{(m)}$ with the same accuracy (\Cref{fig:test1TrueMean}), showing that the UQ-VAE approach is effective in estimating parameters in a noisy setting. The relative errors on the true diffusion coefficients are computed by substituting $\E_\mathrm{p_\mathrm{U|Y}(\mathbf u|\mathbf y^\mathrm{(m)})}[\mathbf u]$ with $\mathbf u^\mathrm{(m)}$ in \eqref{eq:relErrVec}.
\begin{table}[t!]
\scriptsize
\centering
	\begin{tabular}{|c|c|c|}
          \hline
	& Offline phase & Online phase ($100$ samples)\\
          \hline
	UQ-VAE &  $5.5h$ & $10.74s$\\
	Sobol' sequences &    $0.6h$ & $8880s$\\
	\hline 
	\end{tabular}
\caption{Comparison between the time required to solve Bayesian inverse problems with UQ-VAE or approximating \eqref{eq:trueMean} and \eqref{eq:trueCov} with Sobol' sequences.}
\label{table:times}
\end{table}
\begin{figure}[t!]
	\centering
	\includegraphics[width=\linewidth]{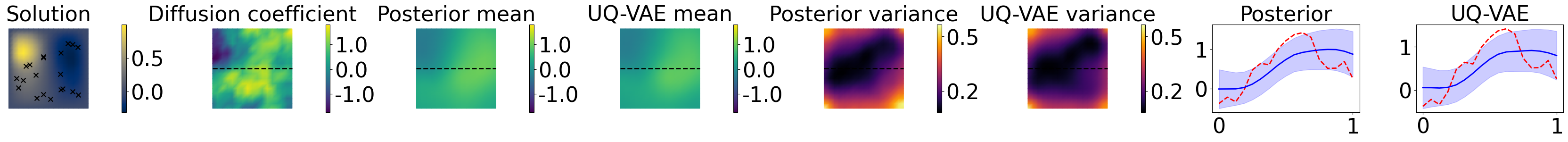}\\
	\includegraphics[width=\linewidth]{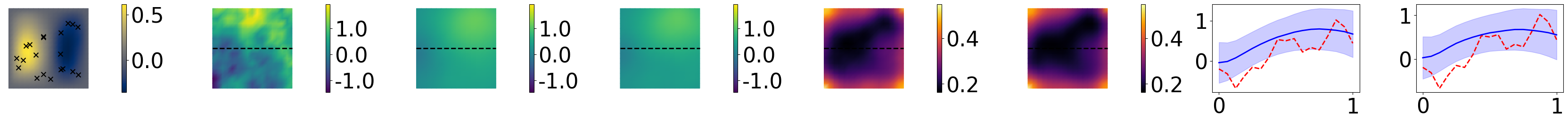} \\
	\includegraphics[width=\linewidth]{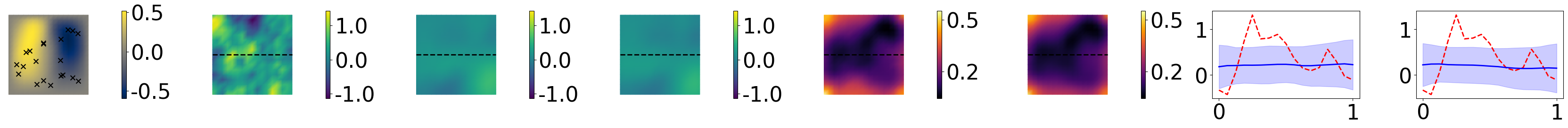}
\caption{True solutions $y^\mathrm{(m)}(\mathbf x)$ with marked observation points, true diffusion coefficients $\mathbf u^\mathrm{(m)}$, posterior and UQ-VAE mean ($\E_\mathrm{p_\mathrm{U|Y}(\mathbf u|\mathbf y^\mathrm{(m)})}[\mathbf u]$ and $\boldsymbol \mu_\mathrm{post}^\mathrm{(m)}$) and variance ($\mathrm {Var}_\mathrm{p_\mathrm{U|Y}(\mathbf u|\mathbf y^\mathrm{(m)})}[\mathbf u]$ and $\Gamma_\mathrm{post}^\mathrm{(m)}$) estimates for three samples. The last two columns are the posterior or UQ-VAE means (continue blue lines) and true diffusion coefficients (dashed red lines) evaluated at $y = 0.5$. The width of the shaded areas is equal to two times the standard deviation.}\label{fig:test1sample}
\end{figure}

The offline computational time for the UQ-VAE (generating the dataset and training the decoder and the encoder) is $5.5h$, whereas the one for generating $2^\mathrm{19}$ samples to estimate \eqref{eq:trueMean} and \eqref{eq:trueCov} is $0.6h$ (\Cref{table:times}). The UQ-VAE approach is advantageous on the online phase, when the Bayesian inverse problems are solved, with an average elapsed time of $0.1s$, compared to the $88.8s$ required to estimate \eqref{eq:trueMean} and \eqref{eq:trueCov}.
\begin{figure}[t!]
	\centering
 	\includegraphics[width=\linewidth, keepaspectratio]{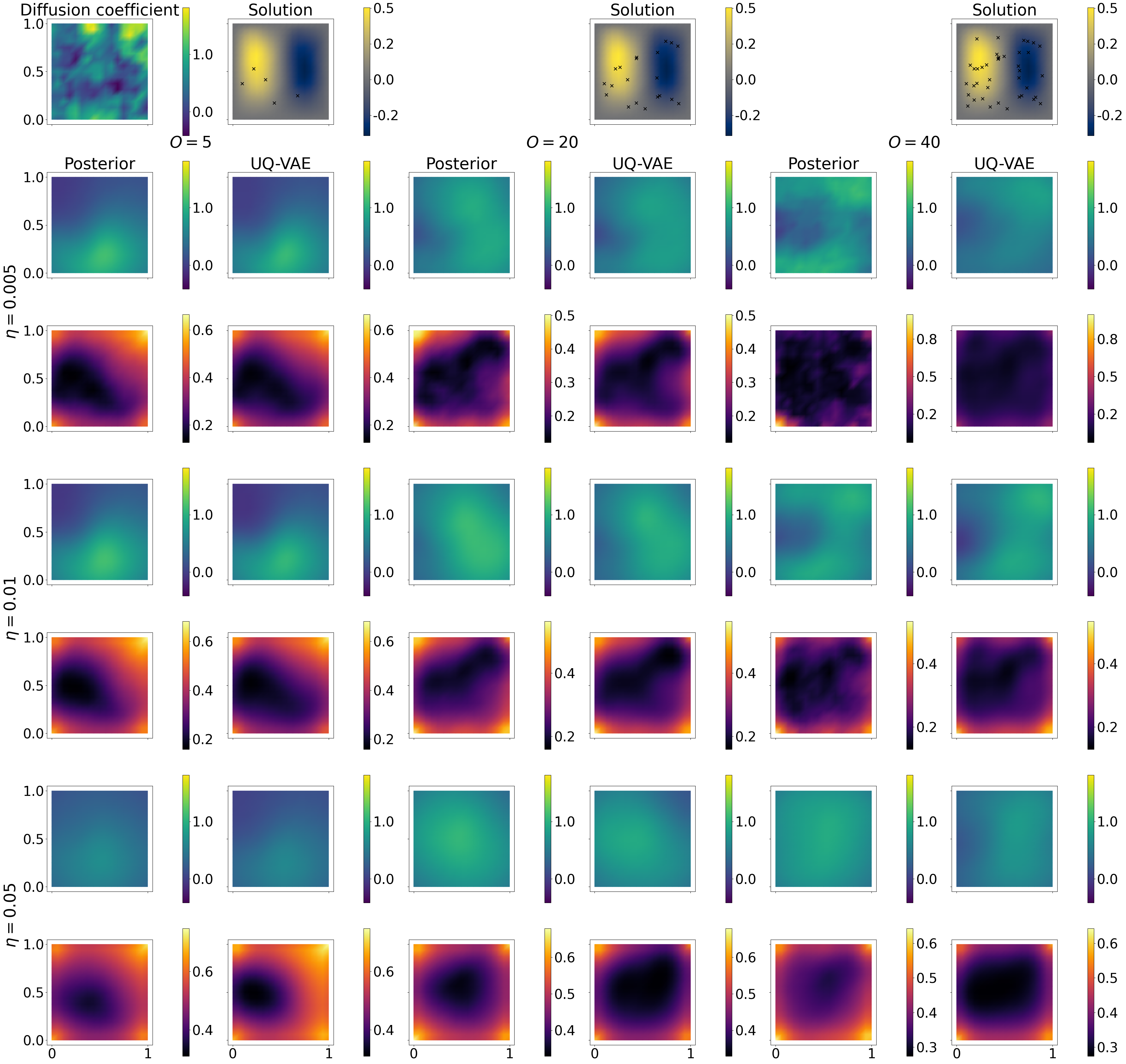}
	\caption{Posterior and UQ-VAE estimates for a test sample, with noise level $\eta$ varying across rows and number of observation points $O$ varying across columns. The first row shows the true diffusion coefficient and the corresponding solution with marked observation points. Estimated means are shown in the second, fourth and sixth rows, while the third, fifth, seventh rows display the corresponding estimated variances.}
	\label{fig:observationNoise}
\end{figure}

We analyze, for three test samples, the posterior and UQ-VAE estimates for the mean and variance (\Cref{fig:test1sample}). The three diffusion coefficients are highly irregular and, when the data are not informative enough, the estimates are poor, like in the third case. For the first two cases, both the posterior and UQ-VAE mean are capable of capturing the trend and the magnitude of the true diffusion coefficient. The two estimates are smoother than the true diffusion coefficient. The UQ-VAE variance is comparable to the posterior one. The estimated mean and variance together are able to catch the irregularity of the true diffusion coefficient, both  for the posterior and the UQ-VAE estimates.

We test the UQ-VAE approach on different values for the noise level $\eta$ and the number of observation points $O$ (\Cref{fig:observationNoise}). For all the values of $\eta$ and $O$ the posterior and UQ-VAE mean estimates are in agreement. For $\eta = 0.05$ the UQ-VAE variance estimates slightly differ from the posterior ones, whereas for $\eta = 0.01$ the two variances are similar. This also happens for $\eta = 0.005$ with $O = 5 \lor 20$. In all such cases, the variance is smaller near the observational data. 
For low values of $\eta$ and large number of observations $O$, the posterior distribution is less regular and the variance estimates get worse, even with $2^\mathrm{19}$ Quasi-Monte Carlo samples, making the comparison with the UQ-VAE estimates unfair (as can be seen for $\eta = 0.005$ and $O = 40$). Indeed, as $\eta$ decreases or $O$ increases, the determinant of the noise covariance matrix $\Gamma_\mathrm{E}$ decreases. The reduction of the regularity of $\Gamma_\mathrm{E}$ corresponds to a lower regularity for $p_\mathrm{U|Y}(\mathbf u|\mathbf y)$ \eqref{eq:postProb}, due to its dependence on $\norm{\cdot}_\mathrm{\Gamma_\mathrm{E}^\mathrm{-1}}$ and $|\Gamma_\mathrm{E}|^\mathrm{-1/2}$. Therefore,  \eqref{eq:trueMean} and \eqref{eq:trueCov} require an elevated number of samples to be estimated accurately, resulting in a high (and unaffordable) computational cost.

\begin{figure}[t!]
	\centering
 	\includegraphics[width=\linewidth, keepaspectratio]{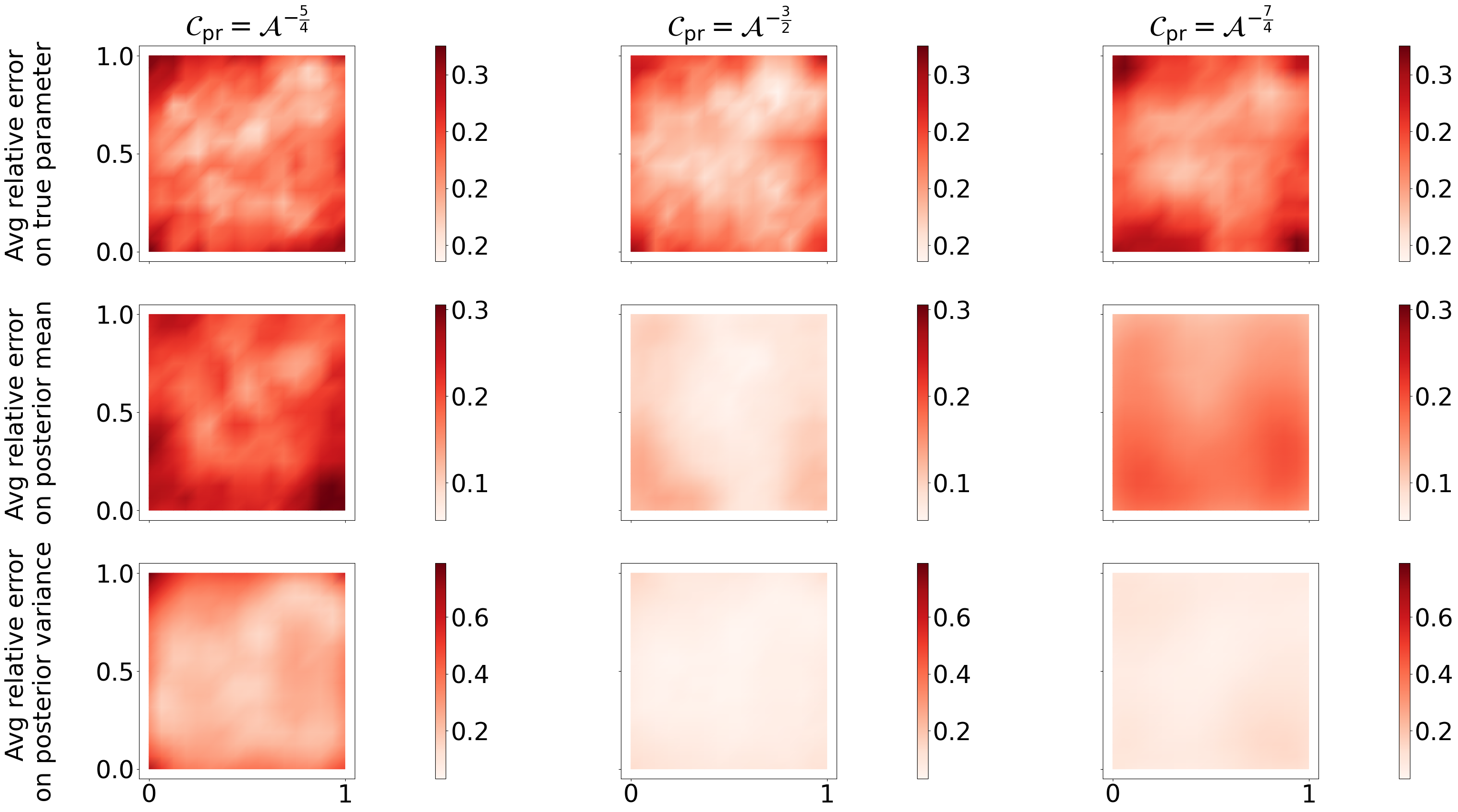}
	\caption{Average relative errors of UQ-VAE estimates on the true parameters and posterior means and variances on $100$ test samples.}
	\label{fig:gammaPrior}
\end{figure}

Finally, we test the UQ-VAE approach for different exponents $\xi$ of the covariance function $\mathcal C_\mathrm{pr}(\mathbf x, \mathbf z)=\mathcal A^\mathrm{-\xi}$ (\Cref{fig:gammaPrior}). We consider three different values of $\xi$ in $(1,2)$ and use the same neural network architecture for all cases. For $\xi=3/2$, we obtain the results already analyzed in \Cref{fig:test1theta,fig:test1TrueMean}. For $\xi=7/4$, the posterior variance is estimated with similar accuracy as in the $\xi=3/2$, while the errors in the estimates of the true parameters and posterior means are slightly higher. This increase in error is due to hyperparameters tuning being performed exclusively for the case $\xi=3/2$. Nonetheless, the errors remain of the same order of magnitude. In the case $\xi=5/4$, the average error on the true parameters is comparable to the other two cases, whereas the errors on the posterior means and variances are higher. As for a low noise level $\eta$ and a high number of observations $O$, smaller values of $\xi$ reduce the regularity of the prior distribution,which in turn decreases the regularity of the posterior distribution. As a result,  \eqref{eq:trueMean} and \eqref{eq:trueCov} would require a large (and computationally infeasible) number of samples to be estimated accurately.

\subsection{Exponential Bayesian inverse problem} \label{sec:nonLinear}
We test the UQ-VAE approach on a non linear Bayesian inverse problem with forward map $\mathcal F(\mathbf u) = e^\mathrm{\mathbf u}$, where the exponentiation is applied element wise. We compare the training time and accuracy of the UQ-VAE trained with loss functions $L_\mathrm{\alpha}$ \eqref{eq:minProblem} and $L_\mathrm{\theta}$ \eqref{eq:minProblem3}. To compute $L_\mathrm{\alpha}$, we use a Sobol' sequence with $K = 2^\mathrm{12}$ samples. In this test, we compare the UQ-VAE estimates with the maximum a posteriori parameter $\mathbf u_\mathrm{MAP}$ \eqref{eq:mapMu} and the Laplace approximation of the covariance matrix $\Gamma_\mathrm{Lap}$ \eqref{eq:Lap} because the true posterior mean \eqref{eq:trueMean} and covariance \eqref{eq:trueCov} are expensive to estimate with an exponential map $\mathcal F$.

We set the dimensions of the parameter and observational data space equal ($D=O$) and vary them from $25$ and $150$. We set $(\mathbf \mu_\mathrm{pr})_\mathrm{i}=-3$, for $i = 1, \dots, D$, and $\Gamma_\mathrm{pr} = 4 I$, where $I \in \R^\mathrm{D \times D}$ is the identity matrix.

We generate a dataset $\{(\mathbf u^\mathrm{(m)}, \mathbf y^\mathrm{(m)}) \}_\mathrm{m=1}^\mathrm{M}$, where $M\in \N$ and $\mathbf u^\mathrm{(m)}$ is sampled from $\mathcal N (\boldsymbol \mu_\mathrm{pr},\Gamma_\mathrm{pr})$. The noiseless observations are denoted by $\tilde{\mathbf y}^\mathrm{(m)} = \mathcal F(\mathbf u^\mathrm{(m)})$. The noisy observations are then defined as $\mathbf y^\mathrm{(m)} = \tilde{\mathbf y}^\mathrm{(m)}+\boldsymbol \epsilon^\mathrm{(m)}$, for $m = 1,\cdots,M$, where $\boldsymbol \epsilon^\mathrm{(m)}$ is sampled from $\mathcal N(\boldsymbol \mu_\mathrm{E},\Gamma_\mathrm{E})$ with $\boldsymbol \mu_\mathrm{E} = \mathbf 0$ and $\Gamma_\mathrm{E}$ diagonal has entries $(\Gamma_\mathrm{E})_\mathrm{i,i} = (\eta \, \underset{m}{\max}\, \abs{\tilde y_\mathrm{i}^\mathrm{(m)}})^\mathrm{2}$, for $i = 1,\cdots, O$ and $\eta \in \R$. We fix $\eta = 0.01$.

We use the true map $\mathcal F$ as the decoder. The encoder consists of $3$ hidden layers with $1000$ neurons each and is trained on a dataset of $100$ samples.
\begin{figure}[t!]
\centering
\includegraphics[width=\linewidth]{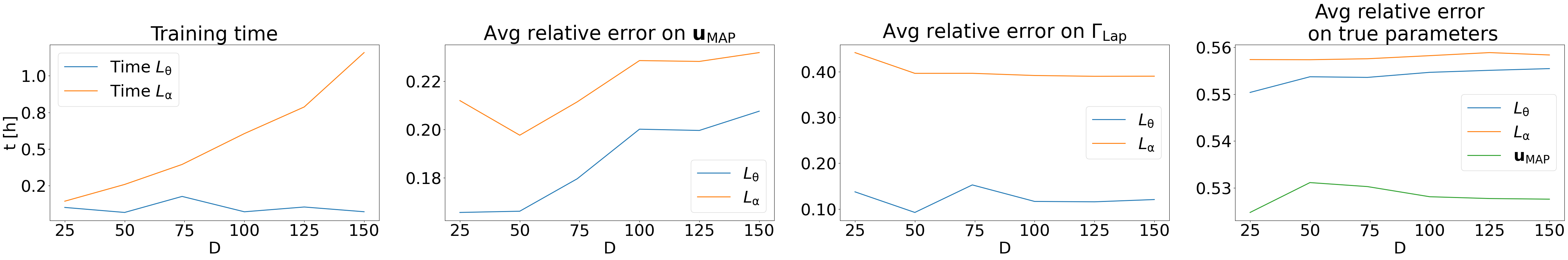}
\caption{Performance comparison in terms of computational time and accuracy of UQ-VAEs (tested on $100$ samples) trained with the loss function $L_\mathrm{\alpha}$ and $L_\mathrm{\theta}$ on the non linear Bayesian inverse problem with forward map $\mathcal F(\mathbf u) = e^\mathrm{\mathbf u}$.}
\label{fig:nonLinear}
\end{figure}%

While the training time of the UQ-VAE using $L_\mathrm{\alpha}$ loss increases with $D$ and $O$, the training time with $L_\mathrm{\theta}$ loss is nearly constant (\Cref{fig:nonLinear}). This is because the computational cost of the training with $L_\mathrm{\theta}$ is mainly driven by the network dimension. In contrast, with $L_\mathrm{\alpha}$, the cost of sampling significantly contributes to the training time, even for small values of $D$.
The average errors (on a dataset of $100$ samples) in the estimation of $\mathbf u_\mathrm{MAP}$, $\Gamma_\mathrm{Lap}$ and true parameters follow a similar trend for both loss functions. Nonetheless, training with $L_\mathrm{\theta}$ yields lower errors across all metrics. The errors in estimating the true parameters with both loss functions are slightly higher than those of $\mathbf u_\mathrm{MAP}$, but still comparable.
\begin{figure}[t!]
\centering
\includegraphics[width=0.5\linewidth]{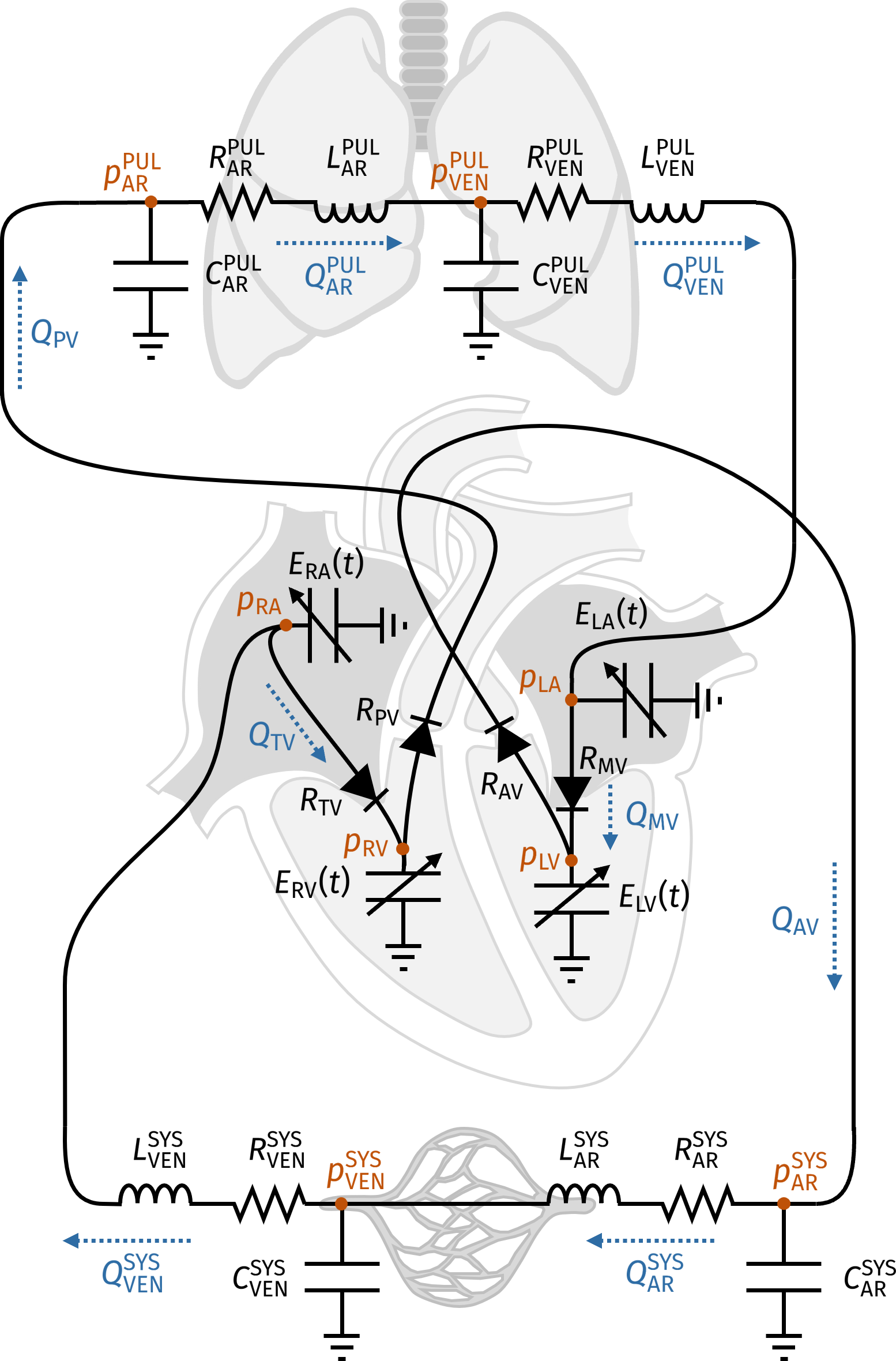}
\caption{0D cardiocirculatory model. We depict pressures and flow rates in red and blue, respectively, and parameters in black.}
\label{fig:lumpedparametercardiovascular}
\end{figure}%

\subsection{0D cardiocirculatory model} \label{sec:cardioTheory}
0D (or lumped-parameter) cardiocirculatory models represent the human cardiovascular system as an electrical circuit: the current corresponds to the blood flow through vessels and valves, the electric potential represents the blood pressure, the electric resistance plays the role of the resistance to blood flow, the capacitance represents the vessel compliance and the inductance corresponds to the blood inertia. In this context, Bayesian inverse problems are used to personalize the model for individual patients and provide insights into their clinical condition, while accounting for the uncertainty in the estimated parameters due to measurement errors (and poor observability of the parameters themselves).

A 0D cardiocirculatory model partitions the cardiovascular system into distinct compartments (e.g. right atrium, systemic arteries/veins). The 0D cardiocirculatory model used here is the one proposed in Ref.~\citen{regazzoni2022cardiac} (\Cref{fig:lumpedparametercardiovascular}). The model consists in the four cardiac chambers, the systemic and pulmonary circulation, split into arterial and venous compartments. \\
The 0D cardiocirculatory model depends on the heart rate (HR), which determines the heartbeat period T$_{\text{HB}}$$=60/$HR, and on the parameters $\check{\mathbf u}$ listed in \Cref{table:params}. A subset $\mathbf u$ of these parameters is used to define the model \eqref{eq:model}, depending on the specific application.\\
\begin{table}[t!]
\begin{center}
\footnotesize
\begin{tabular}{|c|c|c|c|c|}
\hline
& Parameter & Unit & Reference value & Description\\
\hline
\multirow{27}{*}{GSA}&$EA_{\mathrm{LA}}$ & mmHg/mL & $0.2273$ & Left atrial active elastance\\
&$EB_{\mathrm{LA}}$ & mmHg/mL & $0.209$ & Left atrial passive elastance\\
&$V_{\mathrm{U,LA}}$ & mL & $2.0$ & Left atrial unstressed volume\\
&$EA_{\mathrm{LV}}$ & mmHg/mL & $3.0391$ & Left ventricular active elastance\\
&$EB_{\mathrm{LV}}$ & mmHg/mL & $0.10$ & Left ventricular passive elastance\\
&$V_{\mathrm{U,LV}}$ & mL & $16.0$ & Left ventricular unstressed volume\\
&$EA_{\mathrm{RA}}$ & mmHg/mL & $0.0429$ & Right atrial active elastance\\
&$EB_{\mathrm{RA}}$ & mmHg/mL & $0.0636$ & Right atrial passive elastance\\
&$V_{\mathrm{U,RA}}$ & mL & $2.0$ & Right atrial unstressed volume\\
&$EA_{\mathrm{RV}}$ & mmHg/mL & $0.6683$ & Right ventricular active elastance\\
&$EB_{\mathrm{RV}}$ & mmHg/mL & $0.07$ & Right ventricular passive elastance\\
&$V_{\mathrm{U,RV}}$ & mL & $16.0$ & Right ventricular unstressed volume\\
&$R_{\mathrm{min}}$ & mmHg$\cdot$ s/mL & 0.0075 & Minimal valve resistance\\
&$R_{\mathrm{max}}$ & mmHg$\cdot$ s/mL & 75006.2 & Maximal valve resistance\\
&$R_{\mathrm{AR}}^{\mathrm{SYS}}$ & mmHg$\cdot$ s/mL & 0.588 & Systemic arterial resistance\\
&$C_{\mathrm{AR}}^{\mathrm{SYS}}$ & mL/mmHg & 0.96 & Systemic arterial compliance\\
&$L_{\mathrm{AR}}^{\mathrm{SYS}}$ & mmHg$\cdot$ s$^\mathrm{2}$/mL & 0.005 & Systemic arterial inertia\\
&$R_{\mathrm{VEN}}^{\mathrm{SYS}}$ & mmHg$\cdot$ s/mL & 0.352 & Systemic venous resistance\\
&$C_{\mathrm{VEN}}^{\mathrm{SYS}}$ & mL/mmHg & 60.0 & Systemic venous compliance\\
&$L_{\mathrm{VEN}}^{\mathrm{SYS}}$ & mmHg$\cdot$ s$^\mathrm{2}$/mL & 0.0005 & Systemic venous inertia\\
&$R_{\mathrm{AR}}^{\mathrm{PUL}}$ & mmHg$\cdot$ s/mL & 0.104 & Pulmonary arterial resistance\\
&$C_{\mathrm{AR}}^{\mathrm{PUL}}$ & mL/mmHg &5.0 & Pulmonary arterial compliance\\
&$L_{\mathrm{AR}}^{\mathrm{PUL}}$ & mmHg$\cdot$ s$^\mathrm{2}$/mL & 0.0005 & Pulmonary arterial inertia\\
&$R_{\mathrm{VEN}}^{\mathrm{PUL}}$ & mmHg$\cdot$ s/mL & 0.0105 & Pulmonary venous resistance\\
&$C_{\mathrm{VEN}}^{\mathrm{PUL}}$ & mL/mmHg & 16.0 & Pulmonary venous compliance\\
&$L_{\mathrm{VEN}}^{\mathrm{PUL}}$ & mmHg$\cdot$ s$^\mathrm{2}$/mL & 0.0005 & Pulmonary venous inertia\\
&$HR$ & s & $75$ & Heart rate\\
\hline
\multirow{12}{*}{No GSA}&$tC_{\mathrm{LA}}$ & s & $0.79T$$_{\text{HB}}$ & Time of left atrial contraction\\
&$TC_{\mathrm{LA}}$ & s & $0.11$T$_{\text{HB}}$ & Duration of left atrial contraction\\
&$tR_{\mathrm{LA}}$ & s & $tC_{\mathrm{LA}}+TC_{\mathrm{LA}}$ & Time of left atrial relaxation\\
&$TR_{\mathrm{LA}}$ & s & $0.8$T$_{\text{HB}}$ & Duration of left atrial relaxation\\
&$tC_{\mathrm{LV}}$ & s & $0.0$ & Time of left ventricular contraction\\
&$TC_{\mathrm{LV}}$ & s & $0.35$T$_{\text{HB}}$ & Duration of left ventricular contraction\\
&$tR_{\mathrm{LV}}$ & s & $tC_{\mathrm{LV}}+TC_{\mathrm{LV}}$ & Time of left ventricular relaxation\\
&$TR_{\mathrm{LV}}$ & s & $0.4T$T$_{\text{HB}}$ & Duration of left ventricular relaxation\\
&$tC_{\mathrm{RA}}$ & s & $0.8$T$_{\text{HB}}$ & Time of right atrial contraction\\
&$TC_{\mathrm{RA}}$ & s & $0.1$T$_{\text{HB}}$ & Duration of right atrial contraction\\
&$tR_{\mathrm{RA}}$ & s & $tC_{\mathrm{RA}}+TC_{\mathrm{RA}}$ & Time of left atrial relaxation\\
&$TR_{\mathrm{RA}}$ & s & $0.7$T$_{\text{HB}}$ & Duration of left atrial relaxation\\
&$tC_{\mathrm{RV}}$ & s & $0.0$ & Time of right ventricular contraction\\
&$TC_{\mathrm{RV}}$ & s & $0.3$T$_{\text{HB}}$ & Duration of left ventricular contraction\\
&$tR_{\mathrm{RV}}$ & s & $tC_{\mathrm{RV}}+TC_{\mathrm{RV}}$ & Time of right ventricular relaxation\\
&$TR_{\mathrm{RV}}$ & s & $0.4$T$_{\text{HB}}$ & Duration of left ventricular relaxation\\
\hline
\end{tabular}
\caption{List of parameters and their reference values for an ideal healthy individual}
\label{table:params}
\end{center}
\end{table}%
The hemodynamics of the cardiovascular system are described by the following dynamical system:
\begin{align}
\begin{cases}\label{eq:cardio}
\dot V_\mathrm{LA}(t)=Q_\mathrm{VEN}^\mathrm{PUL}(t)-Q_\mathrm{MV}(t),\\[2pt]
\dot V_\mathrm{LV}(t)=Q_\mathrm{MV}(t)-Q_\mathrm{AV}(t),\\[2pt]
 C_\mathrm{AR}^\mathrm{SYS}\dot p_\mathrm{AR}^\mathrm{SYS}(t)=Q_\mathrm{AV}(t)-Q_\mathrm{AR}^\mathrm{SYS}(t),\\[2pt]
L_\mathrm{AR}^\mathrm{SYS}\dot Q_\mathrm{AR}^\mathrm{SYS}(t)=-R_\mathrm{AR}^\mathrm{SYS}Q_\mathrm{AR}^\mathrm{SYS}(t)+p_\mathrm{AR}^\mathrm{SYS}(t)-p_\mathrm{VEN}^\mathrm{SYS}(t),\\[2pt]
C_\mathrm{VEN}^\mathrm{SYS}\dot p_\mathrm{VEN}^\mathrm{SYS}(t)=Q_\mathrm{AR}^\mathrm{SYS}(t)-Q_\mathrm{VEN}^\mathrm{SYS}(t),\\[2pt]
L_\mathrm{VEN}^\mathrm{SYS} \dot Q_\mathrm{VEN}^\mathrm{SYS}(t)=-R_\mathrm{VEN}^\mathrm{SYS}Q_\mathrm{VEN}^\mathrm{SYS}(t)+p_\mathrm{VEN}^\mathrm{SYS}(t)-p_\mathrm{RA}(t),\\[2pt]
\dot V_\mathrm{RA}(t)=Q_\mathrm{VEN}^\mathrm{SYS}(t)-Q_\mathrm{TV}(t),\\[2pt]
\dot V_\mathrm{RV}(t)=Q_\mathrm{TV}(t)-Q_\mathrm{PV}(t),\\[2pt]
C_\mathrm{AR}^\mathrm{PUL}\dot p_\mathrm{AR}^\mathrm{PUL}(t)=Q_\mathrm{PV}(t)-Q_\mathrm{AR}^\mathrm{PUL}(t),\\[2pt]
L_\mathrm{AR}^\mathrm{PUL}\dot Q_\mathrm{AR}^\mathrm{PUL}(t)=-R_\mathrm{AR}^\mathrm{PUL}Q_\mathrm{AR}^\mathrm{PUL}(t)+p_\mathrm{AR}^\mathrm{PUL}(t)-p_\mathrm{VEN}^\mathrm{PUL}(t),\\[2pt]
C_\mathrm{VEN}^\mathrm{PUL}\dot p_\mathrm{VEN}^\mathrm{PUL}(t)=Q_\mathrm{AR}^\mathrm{PUL}(t)-Q_\mathrm{VEN}^\mathrm{PUL}(t),\\[2pt]
L_\mathrm{VEN}^\mathrm{PUL}\dot Q_\mathrm{VEN}^\mathrm{PUL}(t)=-R_\mathrm{VEN}^\mathrm{PUL}Q_\mathrm{VEN}^\mathrm{PUL}(t)+p_\mathrm{VEN}^\mathrm{PUL}(t)-p_\mathrm{LA}(t),
\end{cases}
\end{align}
coupled with suitable initial conditions $\mathbf x_\mathrm{0}$ and with $t \in (0,T]$. The model unknowns are the volumes of the left atrium ($V_\mathrm{LA}$) and ventricle ($V_\mathrm{LV}$) and of the right atrium ($V_\mathrm{RA}$) and ventricle ($V_\mathrm{RV}$), the circulatory pressures of the systemic arteries ($p_\mathrm{AR}^\mathrm{SYS}$) and veins ($p_\mathrm{VEN}^\mathrm{SYS}$) and of the pulmonary arteries ($p_\mathrm{AR}^\mathrm{PUL}$) and veins ($p_\mathrm{VEN}^\mathrm{PUL}$), the circulatory fluxes of the systemic arteries ($Q_\mathrm{AR}^\mathrm{SYS}$) and veins ($Q_\mathrm{VEN}^\mathrm{SYS}$) and of the pulmonary arteries ($Q_\mathrm{AR}^\mathrm{PUL}$) and veins ($Q_\mathrm{VEN}^\mathrm{PUL}$).\\
The fluxes through the valves (tricuspid, pulmonary, mitral and aortic) depend on the pressure jump between the upstream and downstream compartments and the resistance through the valves is given by the non linear function $R_\mathrm{valve}$ as follows:
\begin{gather}
Q_{\mathrm{TV}}(t) = Q_{\mathrm{valve}}(p_{\mathrm{RA}}(t)-p_{\mathrm{RV}}(t)), \quad Q_{\mathrm{MV}}(t) = Q_{\mathrm{valve}}(p_{\mathrm{LA}}(t)-p_{\mathrm{LV}}(t)),\\
Q_{\mathrm{PV}}(t) = Q_{\mathrm{valve}}(p_{\mathrm{RV}}(t)-p_{\mathrm{AR}}^{\mathrm{PUL}}(t)), \quad Q_{\mathrm{AV}}(t) = Q_{\mathrm{valve}}(p_{\mathrm{LV}}(t)-p_{\mathrm{AR}}^{\mathrm{SYS}}(t)),\\
Q_{\mathrm{valve}}(\Delta p) = \frac{\Delta p}{R_{\mathrm{valve}}(\Delta p)},\quad
R_{\mathrm{valve}}(\Delta p) = \sqrt{R_\mathrm{max}R_\mathrm{min}}\left(\frac{R_\mathrm{max}}{R_\mathrm{min}}\right)^\mathrm{\frac{\arctan \left(-100\pi \Delta p\right)}{\pi}}.
\end{gather}
Each cardiac chamber is described as a pressure generator and depends on the periodic function $e_\mathrm{c}(t)$ that models the periodicity of each heartbeat.
\begin{gather}
p_\mathrm{c}(t) = E_\mathrm{c}(t)(V_\mathrm{c}(t)-V_\mathrm{U,c}),\\
E_{\mathrm{c}}(t) = EB_{\mathrm{c}}+EA_{\mathrm{c}}e_{\mathrm{c}}(t),\\
e_{\mathrm{c}}(t) =
\begin{cases}
\frac{1}{2}\left[1-\cos\left(\frac{\pi}{TC_{\mathrm{c}}}\mathrm{mod}\left(t-tC_{\mathrm{c}},T_{\mathrm{HB}} \right) \right)\right] \qquad &\text{if } \mathrm{mod}\left(t-tC_{\mathrm{c}},T_{\mathrm{HB}}\right)<TC_{\mathrm{c}},\\[5pt]
\frac{1}{2}\left[1+\cos\left(\frac{\pi}{TR_{\mathrm{c}}}\mathrm{mod}\left(t-tR_{\mathrm{c}},T_{\mathrm{HB}} \right) \right)\right] \qquad &\text{if } \mathrm{mod}\left(t-tR_{\mathrm{c}},T_{\mathrm{HB}}\right)<TR_{\mathrm{c}},\\[5pt]
0 & \text{otherwise},
\end{cases}
\end{gather}
for $c \in \{LA,LV,RA,RV\}$.

We solve the discrete dynamical system for $25$ heartbeats ($T = 25T_\mathrm{HB}$) to approach its periodic state and focus on the last heartbeat. We compute the model outputs $\check{\mathbf y}$ as functions (e.g. the maximum or the mean) of the unknowns $\mathbf x$ and the parameters $\check{\mathbf u}$:
\begin{align}
\check{\mathbf y} = \mathbf g(\{ \mathbf x(t;\mathbf x_\mathrm{0}, \check{\mathbf u}), t \in [24T_\mathrm{HB},25T_\mathrm{HB}]\};\check{\mathbf u}).
\end{align}
The model outputs are defined in \Cref{table:qoi}, where a subscript \say{I-} refers to the indexed volumes, i.e., normalized by the body surface area (BSA). In general, a subset $\mathbf y$ of the model outputs is used as observational data in model \eqref{eq:model}, depending on the specific application.
\begin{table}[t!]
\scriptsize
\makebox[\textwidth][c]{
	\begin{tabular}{|c|c|c|c|c|}
	\hline
	Model output & Unit & Range & Model value & Description\\
	\hline
	$LA_{\mathrm{I-Vmax}}$ & mL/m$^2$ & [16,34]\cite{lang2015recommendations} & 33.1 & Indexed maximal left atrial volume\\
	$LA_\mathrm{Pmax}$ & mmHg & [6,20]\cite{hurst1990heart} & 12.0 & Maximal left atrial pressure\\
	$LA_\mathrm{Pmin}$ & mmHg & [-2,9]\cite{hurst1990heart} & 8.1 & Minimal left atrial pressure\\
	$LA_\mathrm{Pmean}$ & mmHg & [4,12]\cite{hurst1990heart} & 10.9 & Mean left atrial pressure\\
	$LV_{\mathrm{SV}}$ & mL & [30,80]\cite{lang2015recommendations} & 67.6 & Left ventricular stroke volume\\
	$CI$ & L/min/m$^2$ & [2.8,4.2]\cite{hurst1990heart} & 2.8 & Cardiac index\\
	$LV_{\mathrm{I-EDV}}$ & mL/m$^2$ & [50,90]\cite{hurst1990heart} & 66.0 & Indexed left ventricular end diastolic volume\\
	$LV_{\mathrm{ESV}}$ &  mL & [18,52]\cite{lang2015recommendations} & 50.6 & Left ventricular end systolic volume\\
	$LV_{\mathrm{EF}}$ & \% & [53,73]\cite{lang2015recommendations} & 57.2 & Left ventricular ejection fraction\\
	$LV_{\mathrm{Pmax}}$ & mmHg & [90,140]\cite{hurst1990heart} & 113.2 & Maximal left ventricular pressure\\
	$LV_{\mathrm{Pmin}}$ & mmHg & [4,12]\cite{hurst1990heart} & 5.6 & Minimal left ventricular pressure\\
	$RA_\mathrm{Pmax}$ & mmHg & [2,14]\cite{hurst1990heart} & 8.4 & Maximal right atrial pressure\\
	$RA_\mathrm{Pmin}$ & mmHg & [-2,6]\cite{hurst1990heart} & 5.7 & Minimal right atrial pressure\\
	$RA_\mathrm{Pmean}$ & mmHg & [-1,8]\cite{hurst1990heart} & 7.2 & Mean right atrial pressure\\
	$RV_{\mathrm{I-EDV}}$ &  mL/m$^2$ & [44,80]\cite{rudski2010guidelines} & 65.8 & Indexed right ventricular end diastolic volume\\
	$RV_{\mathrm{I-ESV}}$ &  mL/m$^2$ & [19,46]\cite{rudski2010guidelines} & 28.1 & Indexed right ventricular end systolic volume\\
	$RV_{\mathrm{EF}}$ & \% & [44,71]\cite{rudski2010guidelines} & 57.3 & Right ventricular ejection fraction\\
	$RV_{\mathrm{Pmax}}$ & mmHg & [15,28]\cite{hurst1990heart} & 27.4 & Maximal right ventricular pressure\\
	$RV_{\mathrm{Pmin}}$ & mmHg & [0,8]\cite{hurst1990heart} & 3.9 & Minimal right ventricular pressure\\
	$SAP_{\mathrm{max}}$ &  mmHg & [-,140]\cite{lang2015recommendations} & 112.2 & Systolic systemic arterial pressure\\
	$SAP_{\mathrm{min}}$ &  mmHg & [-,80]\cite{lang2015recommendations} & 62.0 & Diastolic systemic arterial pressure\\
	$PAP_{\mathrm{max}}$ & mmHg & [15,28]\cite{hurst1990heart} & 25.8 & Systolic pulmonary arterial pressure\\
	$PAP_{\mathrm{min}}$ & mmHg & [5,16]\cite{hurst1990heart} & 16.0 & Diastolic pulmonary arterial pressure\\
	$PAP_{\mathrm{mean}}$ & mmHg & [10,22]\cite{hurst1990heart} & 20.6 & Mean pulmonary arterial pressure\\
	$PWP_{\mathrm{min}}$ & mmHg & [1,12]\cite{hurst1990heart} & 11.2 & Minimal pulmonary wedge pressure\\
	$PWP_{\mathrm{mean}}$ & mmHg & [6,15]\cite{hurst1990heart} & 11.8 &Mean pulmonary wedge pressure \\
	$SVR$ & mmHg $\cdot$ min/L & [11.3,17.5]\cite{hurst1990heart} & 15.7 &Systemic vascular resistance\\
	$PVR$ & mmHg $\cdot$ min/L & [1.9,3.1]\cite{hurst1990heart} & 1.9 &Pulmonary vascular resistance \\	
	\hline
	\end{tabular}}

\caption{List of model outputs, their units of measure, the echocardiographic ranges for a healthy individual and the values returned by the 0D cardiocirculatory model with the reference setting of parameters.}
\label{table:qoi}
\end{table}

We determine the reference setting of parameters $\check{\mathbf u}^\mathrm{R}$ to represent an ideal healthy individual (\Cref{table:params}). Specifically, we set HR $=75$ bpm \cite{zhang2009heart} and adjust the values of the other parameters based on modifications of literature values \cite{albanese2016integrated,dede2021modeling,tonini2025two}, ensuring that the model outputs $\check{\mathbf y}$ fall within echocardiograpic ranges characteristic of a healthy individual (\Cref{table:qoi}).

We perform a global sensitivity analysis (GSA) on the $\check{\mathbf u}$ parameters of the 0D cardiocirculatory model indicated in \Cref{table:params} to identify those that most influence a set of observational data $\mathbf y$. This allows to reduce the dimensionality of the parameter space by focusing on the most significant parameters. We compute total Sobol' indices \cite{sobol1993sensitivity} to evaluate the impact of each model parameter on the variance of each observational data $y_\mathrm{j}$, for $j = 1,\dots, O$, accounting for both first-order and higher-order interactions among parameters:
\begin{align}
S_\mathrm{i}^\mathrm{j,T} = 1- \frac{\mathrm {Var}_\mathrm{\check p(\check{\mathbf u}_\mathrm{\sim i})}[\E_\mathrm{\check p(\check u_\mathrm{i})}[y_\mathrm{j}|\check{\mathbf u}_\mathrm{\sim i}]]}{\mathrm {Var}_\mathrm{\check p(\check {\mathbf u})}[y_\mathrm{j}]},
\end{align}
where $\check p$ is the sampling pdf of the parameters and $\check{\mathbf u}_\mathrm{\sim i}$ represents the set of all parameters except for $\check u_\mathrm{i}$. We estimate the total Sobol' indices using Saltelli's method \cite{saltelli2010variance}, allowing the parameters $\check{\mathbf u}$ to vary in a hypercube. The number of samples for the Saltelli's method increases linearly with the number of parameters $N_\mathrm{u}$ and a user-defined value $N\in \N$ as $2N(N_\mathrm{u}+1)$. We set $N = 2^\mathrm{12}$. Since the sensitivity analysis depends on the parameter ranges, we will specify them for the specific application (\cref{sec:test2,sec:test3}). We will solve Bayesian inverse problems for the parameters $\mathbf u \in \R^\mathrm{D}$ associated with at least one Sobol index greater than $0.1$. All the parameters not selected will be kept fixed at their reference value.

\emph{Remark.} The heart rate ($HR$) is usually measured in clinical practice. In our framework there are two possibilities to treat $HR$: ignoring the data and estimating it by means of the UQ-VAE; developing one decoder that can handle different values for $HR$. Even if the second approach reduces the dimension of the parameter space by one, it would require the training of an encoder for each value of $HR$ because the map $\mathcal F$ of the mathematical model \eqref{eq:model} would change with $HR$. Therefore, we restrain ourselves to the first approach.

\subsubsection{Systemic hypertension test case}  \label{sec:test2}
We solve the Bayesian inverse problem for the 0D cardiocirculatory model in the context of systemic hypertension \cite{world2023global}. The chosen observational data are:
\begin{align*}
\mathbf y = [CI, LV_\mathrm{I-ESV}, LV_\mathrm{EF}, LV_\mathrm{Pmax}, SAP_\mathrm{max}, SAP_\mathrm{min}].
\end{align*}
These observational data are directly affected by pressure increases in the systemic circulation.

\begin{figure}[t!]
	\centering
 	\includegraphics[width=\linewidth, keepaspectratio]{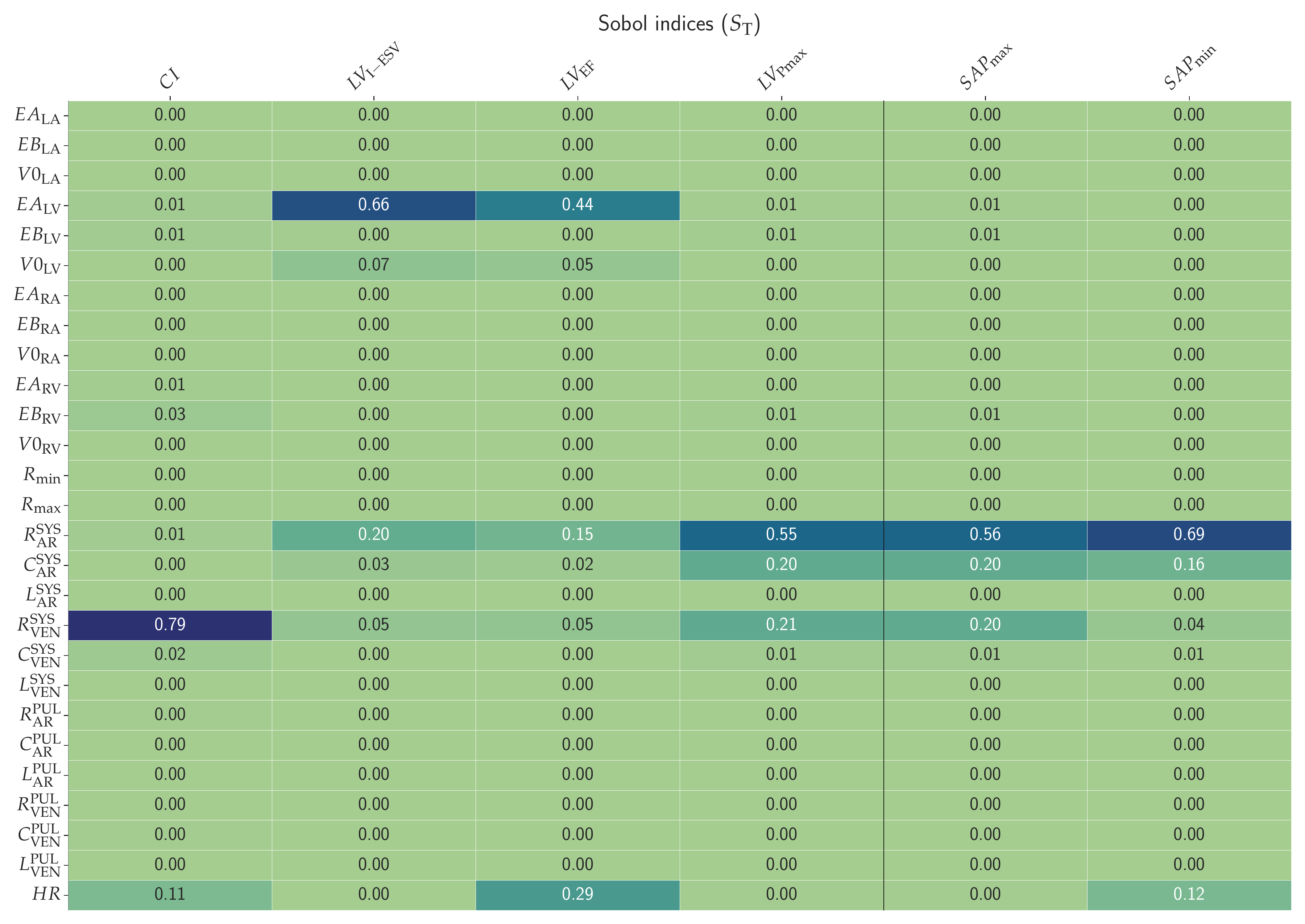}
	\caption{Total Sobol' indices for 0D cardiocirculatory model in the hypertension case.}
	\label{fig:GSAhyper}
\end{figure}

The hypercube used for the sensitivity analysis is built around the reference parameter setting $\check{\mathbf u}^\mathrm{R}$, with variations of $\pm 25 \%$ relative to $\check{\mathbf u}^\mathrm{R}$ to account for intervariability between individuals. Additionally, the parameters listed in \Cref{table:hyper} allow for further adjustments to the lower and upper bounds of the ranges to model systemic hypertension. Based on the total Sobol' indices (\Cref{fig:GSAhyper}), the selected parameters are:
\begin{align*}
\mathbf u=[EA_\mathrm{LV}, R_\mathrm{AR}^\mathrm{SYS}, C_\mathrm{AR}^\mathrm{SYS},R_\mathrm{VEN}^\mathrm{SYS},HR].
\end{align*}
\begin{table}[t]
\scriptsize
\makebox[\textwidth][c]{\renewcommand{\arraystretch}{1.5}
	\begin{tabular}{|c|c|c|c|c|c|c|c|c|c|}
	\hline
	Parameter & $EA_\mathrm{LV}$ & $V_\mathrm{U,LV}$ & $R_\mathrm{AR}^\mathrm{SYS}$ & $C_\mathrm{AR}^\mathrm{SYS} $ & $R_\mathrm{VEN}^{SYS}$ & $C_\mathrm{VEN}^\mathrm{SYS}$ & $R_\mathrm{AR}^\mathrm{PUL}$ & $C_\mathrm{AR}^\mathrm{PUL}$& $HR$ \\
	\hline
	Modification & $+40\%$ & $-10\%$ & $+50\%$ & $-40\%$ & $+5\%$ & $-5\%$ & $+10\%$ & $-10\%$ & $+10\%$   \\
	\hline 
	\end{tabular}}
\caption{List of further modifications in the parameters ranges for the GSA to model hypertension. The modifications affect the lower bound (-) or the upper bound (+) of the ranges.}
\label{table:hyper}
\end{table}%
Since all entries of the parameter vector $\mathbf u$ of the 0D cardiocirculatory model are positive, we work with their logarithms to formulate an unconstrained inverse problem. We assume the parameter vector to be distributed as a lognormal random variable $\log U \sim \mathcal N(\boldsymbol \mu_\mathrm{pr},\Gamma_\mathrm{pr})$. We set the prior mean as $\boldsymbol \mu_\mathrm{pr} = \log \mathbf u^\mathrm{R}$, where $\mathbf u^\mathrm{R}$ is the vector of reference values for $\mathbf u$ and $\log$ is applied element-wise. We let $\Gamma_\mathrm{pr}$ to be a diagonal matrix. During the sensitivity analysis, for $i=1,\dots,D$, each entry $u_\mathrm{i}$ of the parameter vector varies in the interval $[a_\mathrm{i}u_\mathrm{i}^\mathrm{R},b_\mathrm{i}u_\mathrm{i}^\mathrm{R}]$ with $a_\mathrm{i},b_\mathrm{i} \in \R_\mathrm{+},\, a_\mathrm{i}<b_\mathrm{i}$, thus $\log u_\mathrm{i} \in [\log a_\mathrm{i}+\log u_\mathrm{i}^\mathrm{R},\log b_\mathrm{i}+\log u_\mathrm{i}^\mathrm{R}]$. We set $(\Gamma_\mathrm{pr})_\mathrm{i,i} = (\log b_\mathrm{i}-\log a_\mathrm{i})^\mathrm{2}/12$.

We generate a dataset $\{(\mathbf u^\mathrm{(m)}, \mathbf y^\mathrm{(m)}) \}_\mathrm{m=1}^\mathrm{M}$, where $M\in \N$ and $\mathbf u^\mathrm{(m)}$ is sampled from $\mathcal N (\boldsymbol \mu_\mathrm{pr},\Gamma_\mathrm{pr})$. Let $\tilde{\mathbf y}^\mathrm{(m)}$ be the noiseless observations of the solution of \eqref{eq:cardio} with parameter $\mathbf u^\mathrm{(m)}$. Then, $\mathbf y^\mathrm{(m)} = \tilde{\mathbf y}^\mathrm{(m)}+\boldsymbol \epsilon^\mathrm{(m)}$, for $m = 1,\cdots,M$, where $\boldsymbol \epsilon^\mathrm{(m)}$ is sampled from $\mathcal N(\boldsymbol \mu_\mathrm{E},\Gamma_\mathrm{E})$ with $\boldsymbol \mu_\mathrm{E} = \mathbf 0$ and $\Gamma_\mathrm{E}$ diagonal with $(\Gamma_\mathrm{E})_\mathrm{i,i} = (\eta \, \underset{m}{\max}\, \abs{\tilde y_\mathrm{i}^\mathrm{(m)}})^\mathrm{2}$, for $i = 1,\cdots, O$ and $\eta \in \R$. Observe that differently from \cref{sec:test1}, $(\Gamma_\mathrm{E})_\mathrm{i,i}$ depends on the specific observational data because they are associated with different quantities. We fix $\eta = 0.05$.

We use a decoder with $5$ hidden layers and $250$ neurons per layer trained on a dataset of $M = 4096$ samples. The encoder consists of $3$ hidden layers and $250$ neurons per layer and is trained on a subset of $100$ samples of the decoder dataset.

We generate $2^\mathrm{17}$ samples from the distribution $\mathcal N(\boldsymbol \mu_\mathrm{pr},\Gamma_\mathrm{pr})$ to estimate \eqref{eq:trueMean} and \eqref{eq:trueCov}.

\begin{figure}[t!]
	\centering
 	\includegraphics[width=\linewidth,keepaspectratio]{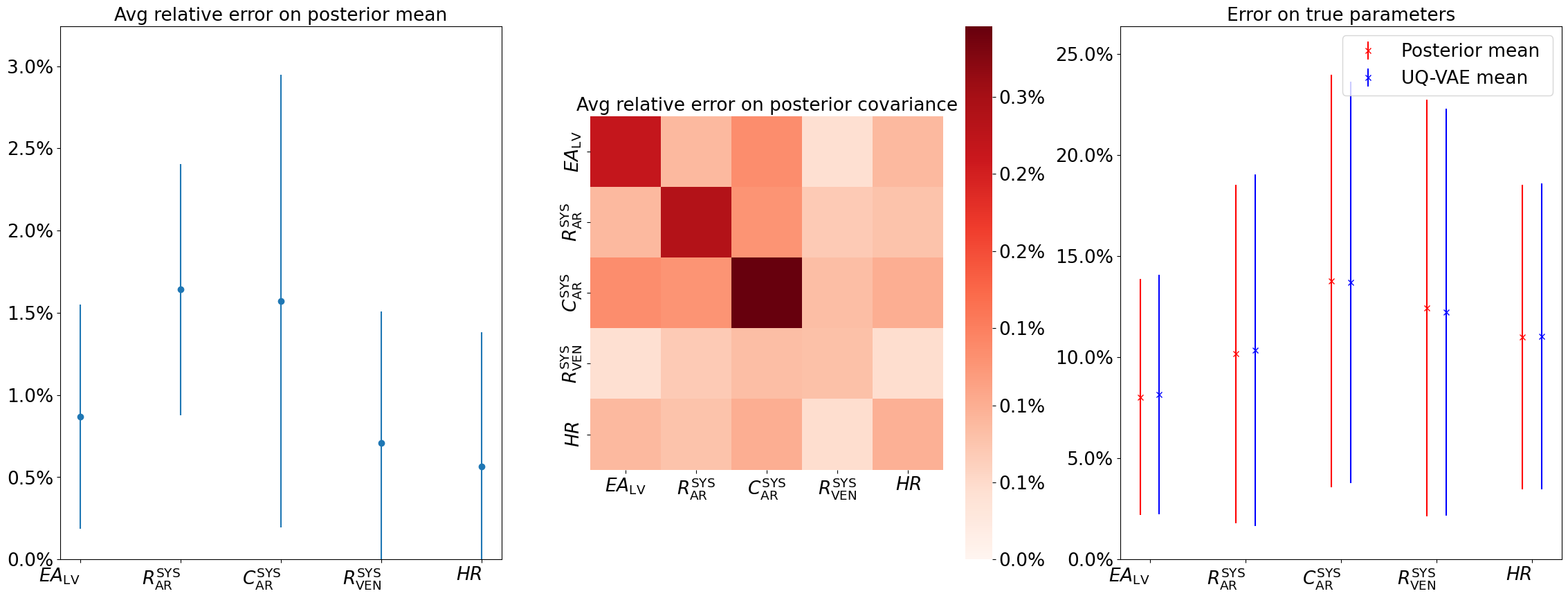}
\caption{From left to right: average and standard deviation of the relative errors of the UQ-VAE mean $\boldsymbol \mu_\mathrm{post}^\mathrm{(m)}$ on the posterior one $\E_\mathrm{p_\mathrm{U|Y}(\mathbf u|\mathbf y^\mathrm{(m)})}[\mathbf u]$; average relative error of the UQ-VAE covariance $\Gamma_\mathrm{post}^\mathrm{(m)}$ on the posterior one $\mathrm {Var}_\mathrm{p_\mathrm{U|Y}(\mathbf u|\mathbf y^\mathrm{(m)})}[\mathbf u]$; average and standard deviations of the relative errors of the posterior and UQ-VAE means on the true parameters $\mathbf u^\mathrm{(m)}$. The test are run for $100$ samples for the systemic hypertension case and the relative errors and standard deviations are computed element-wise.}\label{fig:errorsHypertensionQoI}
\end{figure}

The UQ-VAE mean $\boldsymbol \mu_\mathrm{post}^\mathrm{(m)}$ estimates the posterior mean $\E_\mathrm{p_\mathrm{U|Y}(\mathbf u|\mathbf y^\mathrm{(m)})}[\mathbf u]$ with an average relative error less than $2.0\%$ on $100$ test samples (\Cref{fig:errorsHypertensionQoI}). The UQ-VAE covariance $\Gamma_\mathrm{post}^\mathrm{(m)}$ estimates the posterior covariance $\mathrm {Cov}_\mathrm{p_\mathrm{U|Y}(\mathbf u|\mathbf y^\mathrm{(m)})}(\mathbf u)$ with an average relative error less than $0.5\%$. Both the posterior and the UQ-VAE means achieve similar average relative errors on the true parameters $\mathbf u^\mathrm{(m)}$.

The chosen observational data $\mathbf y$ are related to the practical case of systemic hypertension. We now test the UQ-VAE using $25$ equally spaced observations of all the $12$ unknowns of \eqref{eq:cardio} in $[24T_\mathrm{HB},25T_\mathrm{HB}]$ ($O=300$). We perform a new sensitivity analysis to account for the different observational data and we select the parameters:
\begin{align}
 \mathbf u = [&EA_\mathrm{LA}, EB_\mathrm{LA},EA_\mathrm{LV},EB_\mathrm{LV},EB_\mathrm{RA},EA_\mathrm{RV},EB_\mathrm{RV},V_\mathrm{U,RV},\\
 &R_\mathrm{AR}^\mathrm{SYS}, C_\mathrm{AR}^\mathrm{SYS},R_\mathrm{VEN}^\mathrm{SYS}, C_\mathrm{VEN}^\mathrm{SYS},R_\mathrm{AR}^\mathrm{PUL}, C_\mathrm{AR}^\mathrm{PUL},L_\mathrm{AR}^\mathrm{PUL},HR].
\end{align}
The definition of the prior distribution and of the dataset are analogous to the previous case. The noise distribution has mean $\boldsymbol \mu_\mathrm{E} = \mathbf 0$ and diagonal covariance matrix $\Gamma_\mathrm{E}$. $\Gamma_\mathrm{E}$ is such that all the observations of the same unknown of \eqref{eq:cardio} are related to the same variance obtained by taking the maximum of the absolute values of all the corresponding samples and observations:
\begin{align}
(\Gamma_\mathrm{E})_{i,i} = \left(\eta \, \underset{\substack{m\\ \\1+25\mathrm{mod}(i-1,25)\le k \le 25(\mathrm{mod}(i-1,25)+1)}}{\max}\, \abs{\tilde y_\mathrm{k}^\mathrm{(m)}}\right)^\mathrm{2}. 
\end{align}

We generate $2^\mathrm{19}$ samples from the distribution $\mathcal N(\boldsymbol \mu_\mathrm{pr},\Gamma_\mathrm{pr})$ to estimate \eqref{eq:trueMean} and \eqref{eq:trueCov}.

\begin{figure}[t!]
	\centering
 	\includegraphics[width=0.9\linewidth, keepaspectratio]{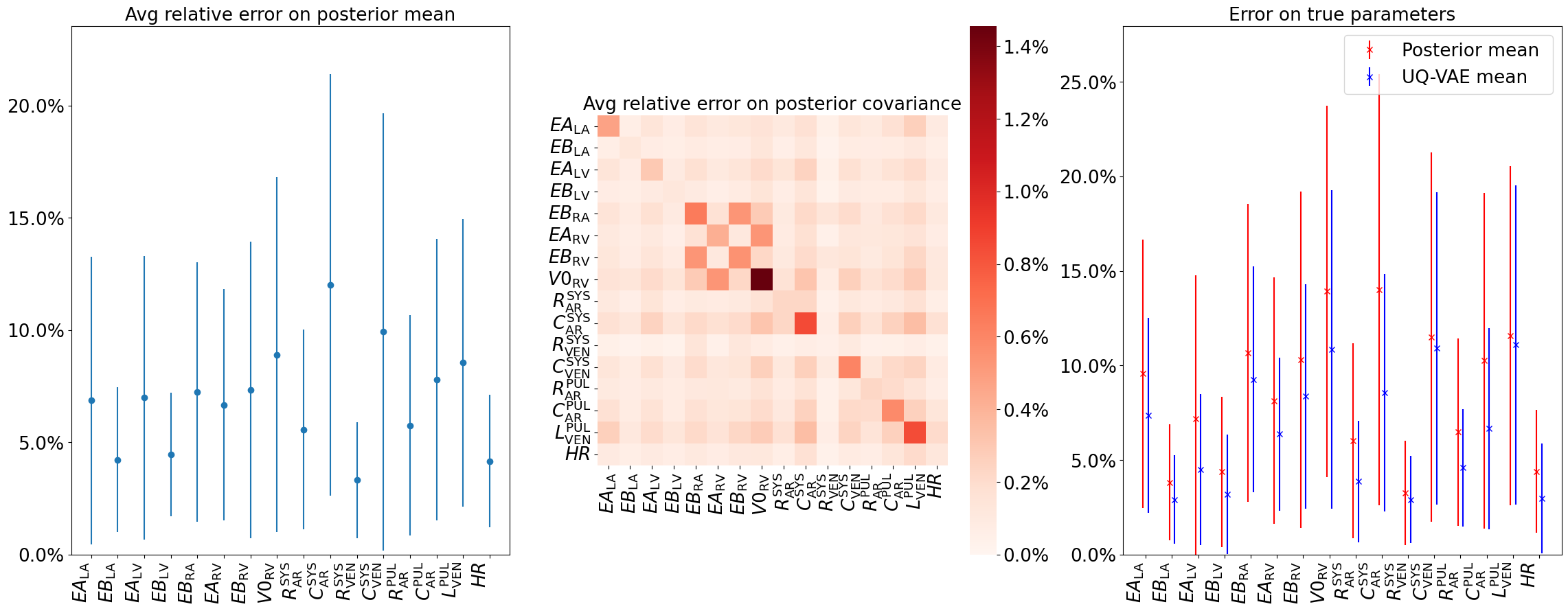}
	\caption{From left to right: average and standard deviation of the relative errors of the UQ-VAE mean $\boldsymbol \mu_\mathrm{post}^\mathrm{(m)}$ on the posterior one $\E_\mathrm{p_\mathrm{U|Y}(\mathbf u|\mathbf y^\mathrm{(m)})}[\mathbf u]$; average relative error of the UQ-VAE covariance $\Gamma_\mathrm{post}^\mathrm{(m)}$ on the posterior one $\mathrm {Var}_\mathrm{p_\mathrm{U|Y}(\mathbf u|\mathbf y^\mathrm{(m)})}[\mathbf u]$; average and standard deviations of the relative errors of the posterior and UQ-VAE means on the true parameters $\mathbf u^\mathrm{(m)}$. The test are run for $100$ samples for the systemic hypertension case and the relative errors and standard deviations are computed element-wise. The data are $25$ equally spaced observations of all the $12$ unknowns of \eqref{eq:cardio}.}
	\label{fig:errorsHypertensionSol}
\end{figure}

In this case, the UQ-VAE mean is more accurate than the posterior mean in estimating the true parameters (\Cref{fig:errorsHypertensionSol}). As in \Cref{fig:observationNoise}, this is because the increase in the number of observational data corresponds to a decrease in the regularity of $p_\mathrm{U|Y}(\mathbf u|\mathbf y)$ and a poorer estimate of \eqref{eq:trueMean} and \eqref{eq:trueCov}. Consequently, also the estimates of the posterior mean and covariance matrix are worse than \Cref{fig:errorsHypertensionQoI} with average relative errors less than $15\%$ and $1.6\%$, respectively. Moreover, the accuracy of the UQ-VAE on the true parameters is unaffected by the higher number of parameters to estimate with respect to the previous case. 

\subsubsection{Ventricular septal defect test case}  \label{sec:test3}
We solve the Bayesian inverse problem for the 0D cardiocirculatory model in the context of ventricular septal defect \cite{spicer2014ventricular,dakkak2024ventricular}. This congenital defect is characterized by a hole in the septum (wall) between the two ventricles. To account for the blood flow through the ventricular septal defect, the 0D cardiocirculatory model requires a modification in the differential equations for the volumes of the left and right ventricles:
\begin{align}\label{eq:mod}
\begin{cases}
\dot V_\mathrm{LV}(t)=Q_\mathrm{MV}(t)-Q_\mathrm{AV}(t)-Q_\mathrm{VSD}(t),\\[2pt]
\dot V_\mathrm{RV}(t)=Q_\mathrm{TV}(t)-Q_\mathrm{PV}(t)+Q_\mathrm{VSD}(t),
\end{cases}
\end{align}
where:
\begin{align}
Q_\mathrm{VSD}(t) = \frac{p_\mathrm{LV}(t)-p_\mathrm{RV}(t)}{R_\mathrm{VSD}},
\end{align}
with $R_\mathrm{VSD}$ representing the resistance to the blood flow through the ventricular septal defect. We use as parameter the radius of the ventricular septal defect $r_\mathrm{VSD}$ with a reference value of $0.6r_\mathrm{AV} = 0.9$ cm, where $r_\mathrm{AV}=1.5$ cm is the aortic annulus radius \cite{dakkak2024ventricular}. We compute $R_\mathrm{VSD}$ by assuming that the length of the ventricular septal defect and the thickness of the aortic annulus are equal to $L$ and applying two times the Poiseuille law:
\begin{align}
R_\mathrm{VSD} = \frac{8 \mu L}{\pi r_\mathrm{VSD}^\mathrm{4}}=R_\mathrm{min}\left(\frac{r_\mathrm{AV}}{r_\mathrm{VSD}}\right)^\mathrm{4},
\end{align}
where $\mu$ denotes the blood viscosity.

The observational data are chosen as follows:
\begin{align*}
\mathbf y =& [LV_\mathrm{Pmax}, LV_\mathrm{Pmin},RA_\mathrm{Pmean},RV_\mathrm{Pmax},RV_\mathrm{Pmin},SAP_\mathrm{max}, SAP_\mathrm{min},\\
&PAP_\mathrm{max},PAP_\mathrm{min}, PAP_\mathrm{mean},PWP_\mathrm{mean}, PVR,Q_\mathrm{P},Q_\mathrm{S}],
\end{align*}
where $Q_\mathrm{P}$ and $Q_\mathrm{S}$ are the mean flows through the pulmonary and aortic valves. These two observational data are particularly important in the case of ventricular septal defects because their ratio $Q_\mathrm{P}/Q_\mathrm{S}$ is typically around $1$ in healthy individuals, indicating that the same amount of blood flows through the pulmonary and systemic circulation. In the case of ventricular septal defect, this ratio increases beyond one, meaning that some of the oxygenated blood from the lungs flows to the right ventricle, bypassing the systemic circulation and failing to supply oxygen to the body's cells.

\begin{figure}[t!]
	\centering
 	\includegraphics[width=0.9\linewidth, keepaspectratio]{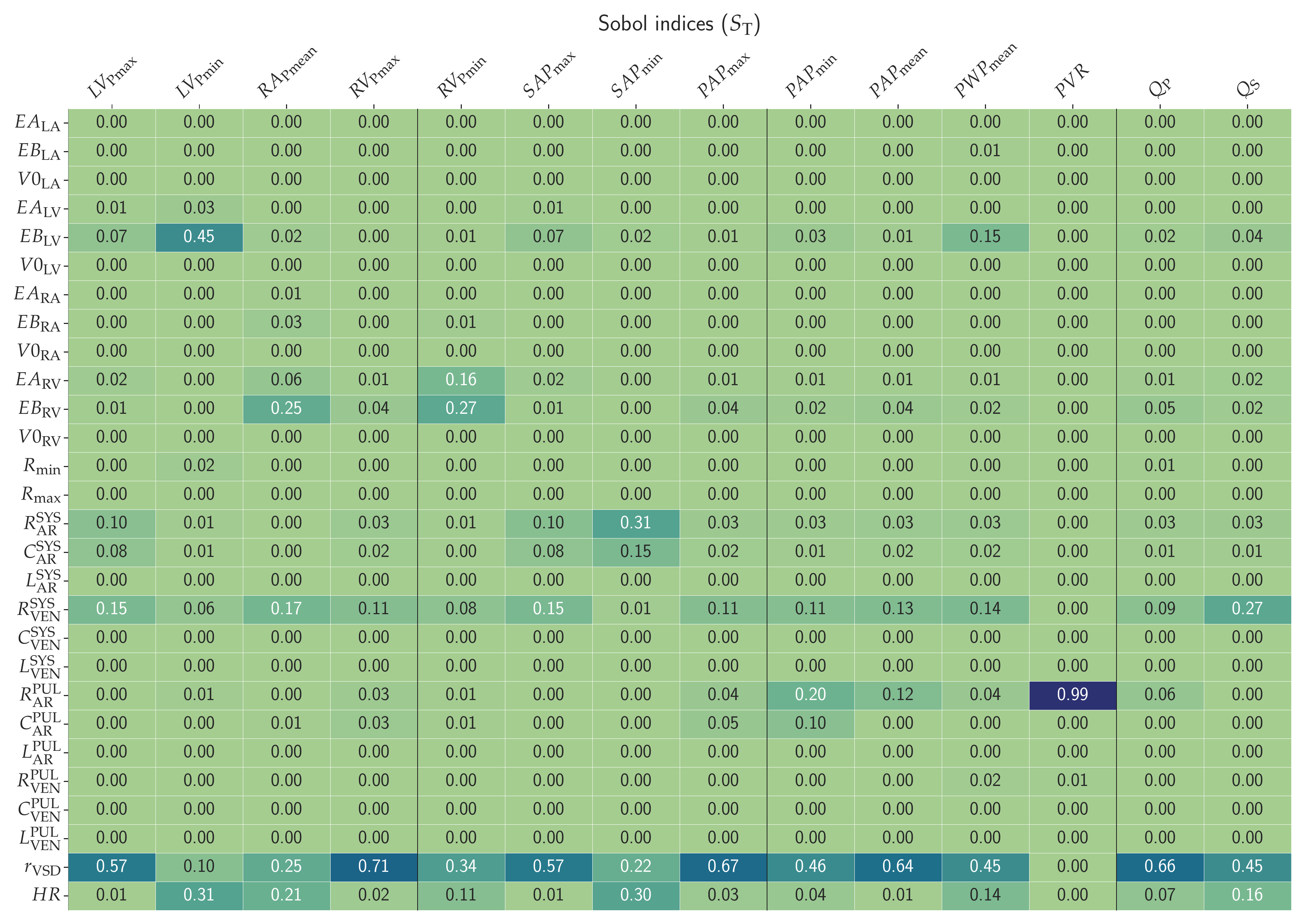}
	\caption{Total Sobol' indices for 0D cardiocirculatory model in the ventricular septal defect case.}
	\label{fig:GSAVSD}
\end{figure}

The hypercube used for the sensitivity analysis is built around the reference setting of parameters $\check{\mathbf u}^\mathrm{R}$, accounting for intervariability between individuals by allowing a variation of $\pm 25 \%$ with respect to $\check{\mathbf u}^\mathrm{R}$. Based on the sensitivity analysis results (\Cref{fig:GSAVSD}), the selected parameters are:
\begin{align*}
 \mathbf u = [EB_\mathrm{LV}, EA_\mathrm{RV}, EB_\mathrm{RV}, R_\mathrm{AR}^\mathrm{SYS}, C_\mathrm{AR}^\mathrm{SYS}, R_\mathrm{VEN}^\mathrm{SYS}, R_\mathrm{AR}^\mathrm{PUL}, r_\mathrm{VSD}, HR].
\end{align*}
Since all entries of the parameter vector $\mathbf u$ of the 0D cardiocirculatory model are positive, we work with their logarithms  to formulate an unconstrained inverse problem. We assume that the parameter vector follows a lognormal distribution $\log U \sim \mathcal N(\boldsymbol \mu_\mathrm{pr},\Gamma_\mathrm{pr})$. We set the prior mean as $\boldsymbol \mu_\mathrm{pr} = \log \mathbf u^\mathrm{R}$, where $\mathbf u^\mathrm{R}$ is the vector of reference values for $\mathbf u$ and $\log$ is applied element-wise. We let $\Gamma_\mathrm{pr}$ to be a diagonal matrix. During the sensitivity analysis, for $i=1,\dots,D$, each entry $u_\mathrm{i}$ of the parameter vector varies within the range $ [0.75u_\mathrm{i}^\mathrm{R},1.25u_\mathrm{i}^\mathrm{R}]$, so that $\log u_\mathrm{i} \in [\log 0.75+\log u_\mathrm{i}^\mathrm{R},\log 1.25+\log u_\mathrm{i}^\mathrm{R}]$. We set $(\Gamma_\mathrm{pr})_\mathrm{i,i} = (\log 1.25-\log 0.75)^\mathrm{2}/12$.

We generate a dataset $\{(\mathbf u^\mathrm{(m)}, \mathbf y^\mathrm{(m)}) \}_\mathrm{m=1}^\mathrm{M}$, where $M\in \N$ and $\mathbf u^\mathrm{(m)}$ is sampled from $\mathcal N (\boldsymbol \mu_\mathrm{pr},\Gamma_\mathrm{pr})$. Let $\tilde{\mathbf y}^\mathrm{(m)}$ be the noiseless observations of the solution of \eqref{eq:cardio} modified with \eqref{eq:mod} with parameter $\mathbf u^\mathrm{(m)}$. Then, $\mathbf y^\mathrm{(m)} = \tilde{\mathbf y}^\mathrm{(m)}+\boldsymbol \epsilon^\mathrm{(m)}$, for $m = 1,\cdots,M$, where $\boldsymbol \epsilon^\mathrm{(m)}$ is sampled from $\mathcal N(\boldsymbol \mu_\mathrm{E},\Gamma_\mathrm{E})$ with $\boldsymbol \mu_\mathrm{E} = \mathbf 0$ and $\Gamma_\mathrm{E}$ diagonal with $(\Gamma_\mathrm{E})_\mathrm{i,i} = (\eta \, \underset{m}{\max}\, \abs {\tilde y_\mathrm{i}^\mathrm{(m)}})^\mathrm{2}$, for $i = 1,\cdots, O$ and $\eta \in \R$. We fix $\eta = 0.05$.

We use a decoder with $5$ hidden layers and $250$ neurons per layer trained on a dataset of $M = 4096$ samples. The encoder consists of $3$ hidden layers and $250$ neurons per layer and is trained on a subset of $100$ samples of the decoder dataset.

We generate $2^\mathrm{17}$ samples from the distribution $\mathcal N(\boldsymbol \mu_\mathrm{pr},\Gamma_\mathrm{pr})$ to estimate \eqref{eq:trueMean} and \eqref{eq:trueCov}.

\begin{figure}[t!]
	\centering
 	\includegraphics[width=0.9\linewidth, keepaspectratio]{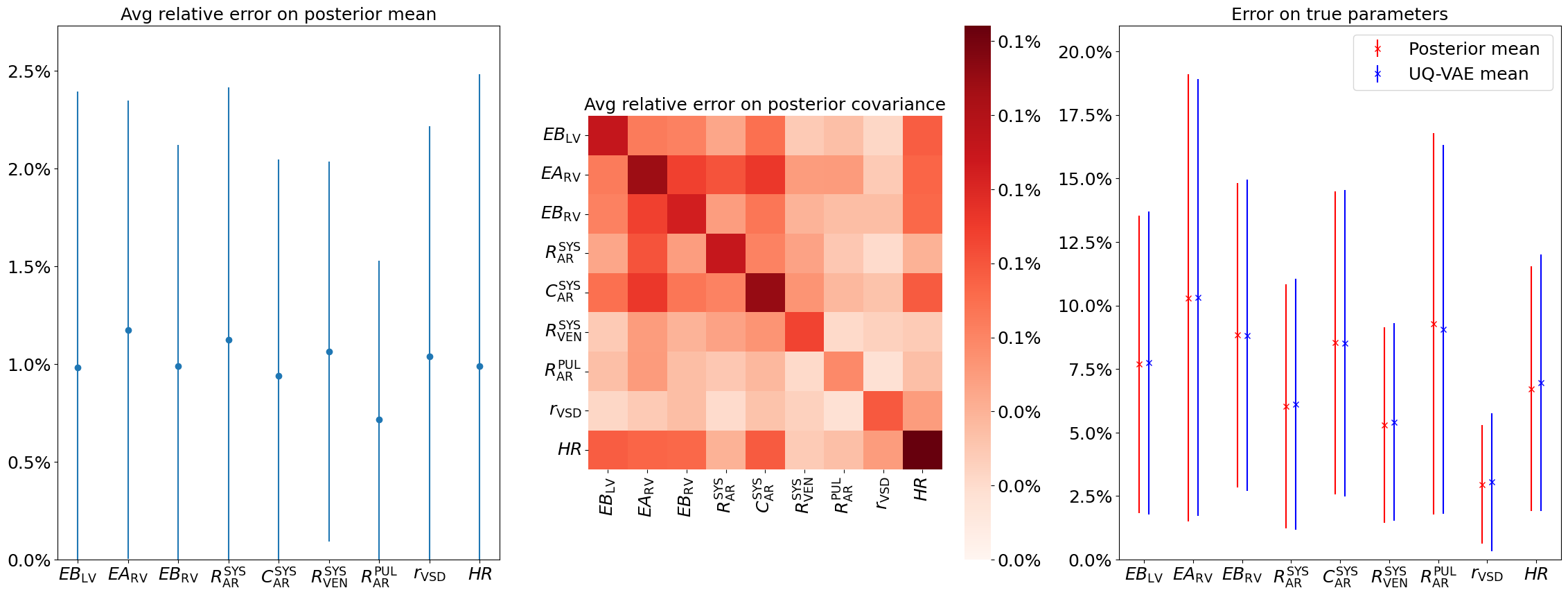}
	\caption{From left to right: average and standard deviation of the relative errors of the UQ-VAE mean $\boldsymbol \mu_\mathrm{post}^\mathrm{(m)}$ on the posterior one $\E_\mathrm{p_\mathrm{U|Y}(\mathbf u|\mathbf y^\mathrm{(m)})}[\mathbf u]$; average relative error of the UQ-VAE covariance $\Gamma_\mathrm{post}^\mathrm{(m)}$ on the posterior one $\mathrm {Var}_\mathrm{p_\mathrm{U|Y}(\mathbf u|\mathbf y^\mathrm{(m)})}[\mathbf u]$; average and standard deviations of the relative errors of the posterior and UQ-VAE means on the true parameters $\mathbf u^\mathrm{(m)}$. The test are run for $100$ samples for the ventricular septal defect case and the relative errors and standard deviations are computed element-wise.}
	\label{fig:errorsVSD}
\end{figure}

The UQ-VAE mean $\boldsymbol \mu_\mathrm{post}^\mathrm{(m)}$ estimates the posterior mean $\E_\mathrm{p_\mathrm{U|Y}(\mathbf u|\mathbf y^\mathrm{(m)})}[\mathbf u]$ with an average relative error less than $2.0\%$ on $100$ test samples (\Cref{fig:errorsVSD}). The UQ-VAE covariance $\Gamma_\mathrm{post}^\mathrm{(m)}$ estimates the posterior covariance $\mathrm {Cov}_\mathrm{p_\mathrm{U|Y}(\mathbf u|\mathbf y^\mathrm{(m)})}(\mathbf u)$ with an average relative error less than $0.2\%$. Both the posterior and the UQ-VAE means achieve similar average relative errors on the true parameters $\mathbf u^\mathrm{(m)}$.

\begin{figure}[t!]
	\centering
 	\includegraphics[width=0.7\linewidth, keepaspectratio]{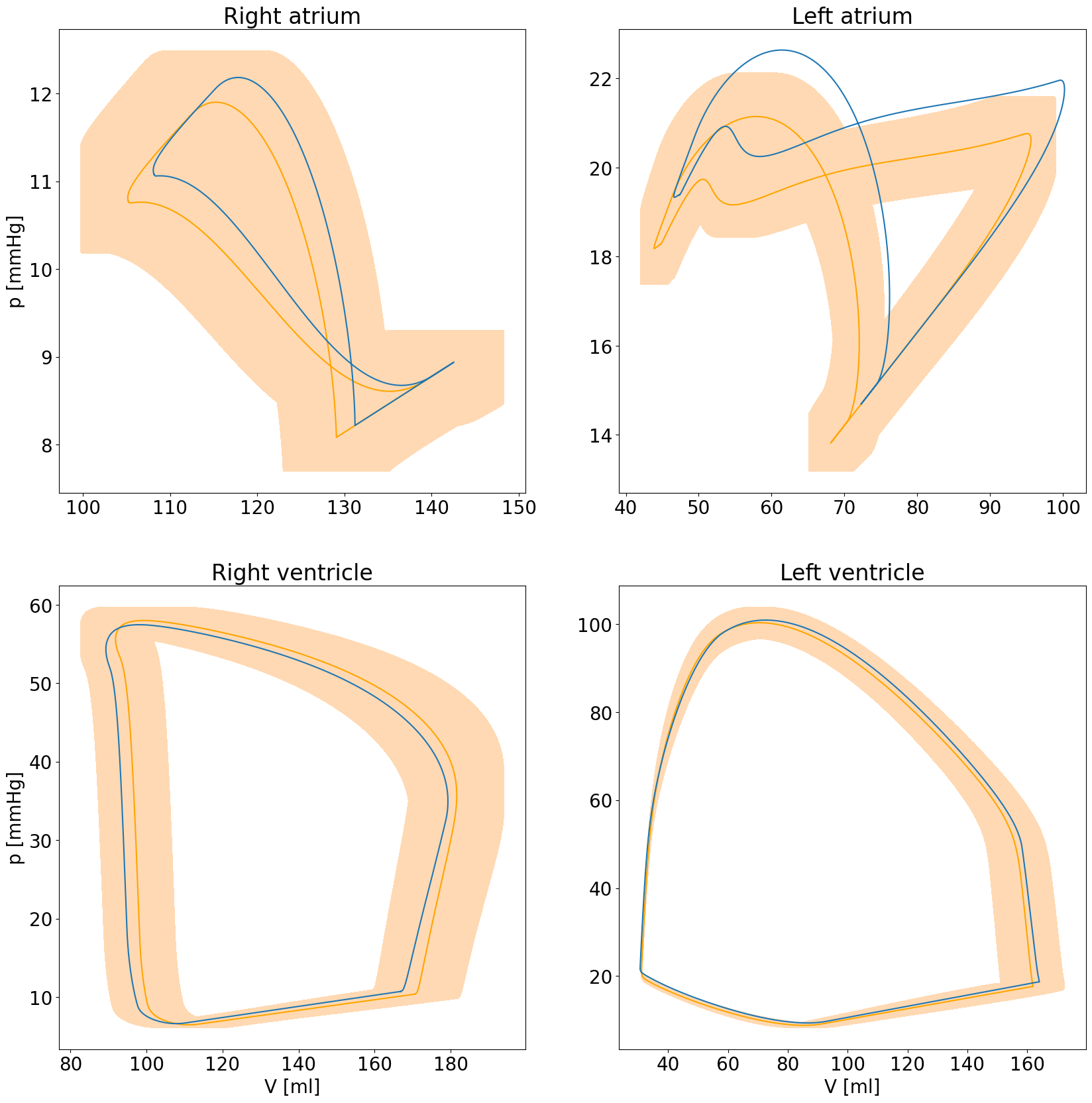}
	\caption{Comparison between the true PV loops (blue) and the mean (orange) of those estimated by sampling the parameters from the approximate posterior distribution $q_\mathrm{\phi}(\mathbf u|\mathbf y)$. The shaded area represents the standard deviation of the estimated PV loops.}
	\label{fig:PVloops}
\end{figure}

We perform forward uncertainty quantification of a single test sample. We use a Sobol' sequence to generate $2^\mathrm{16}$ parameter samples from the approximate posterior distribution $q_\mathrm{\phi}(\mathbf u|\mathbf y) = \mathcal N(\boldsymbol \mu_\mathrm{post},\Gamma_\mathrm{post})$. For each parameter sample, we run the 0D cardiocirculatory model with the VSD and compute the mean and standard deviation of the resulting time transients pressures and volumes in the four cardiac chambers. We compare the means and standard deviation of the estimated pressure-volume (PV) loops with the true PV loops (\Cref{fig:PVloops}). PV loops are indicators of cardiac dysfunction and their accurate estimation provides valuable clinical insights into a patient's condition. The estimated PV loops for the right atrium and both ventricles are highly accurate, with smaller uncertainties observed in the left ventricle. The estimation for the left atrium is less accurate compared to the other chambers. This discrepancy may be attributed to the high estimation error in $EB_\mathrm{LV}$ (\Cref{table:estParam}) which affects both the pressure and volume in the left atrium. Another possible explanation is that the true posterior distribution may not be Gaussian, as we assumed in our approximation, potentially introducing a bias in the estimated distribution.
\begin{table}[t!]
\scriptsize
\centering
\begin{tabular}{|c|c|c|c|c|c|c|c|c|c|}
\hline
& $EB_\mathrm{LV}$ & $EA_\mathrm{RV}$ & $EB_\mathrm{RV}$ & $R_\mathrm{AR}^\mathrm{SYS}$ & $C_\mathrm{AR}^\mathrm{SYS}$ & $R_\mathrm{VEN}^\mathrm{SYS}$ & $R_\mathrm{AR}^\mathrm{PUL}$ & $r_\mathrm{VSD}$ & HR\\
\hline
$\mathbf u$ & $0.11$ &  $0.70$ & $0.087$ & $0.66$ & $1.1$ & $0.32$ & $0.10$ & $0.79$ & $82$\\
$\boldsymbol \mu_\mathrm{post}$ & $0.09$ & $0.64$ & $0.074$ & $0.62$ & $1.0$ & $0.32$ & $0.09$ & $0.82$ & $82$\\
Relative error & $0.14$ & $0.08$ & $0.15$ & $0.07$ & $0.06$ & $0.01$ & $0.10$ & $0.04$ & $0.01$\\
\hline
\end{tabular}
\caption{True and estimated parameter, along with their relative error, for a test sample of the 0D cardiocirculatory model with VSD.}
\label{table:estParam}
\end{table}

\section{Conclusions}\label{sec:concl}
We proposed an improvement (in the form of a novel loss function) of Uncertainty Quantification Variational AutoEncoders (UQ-VAE) for solving Bayesian inverse problems. We introduced a new loss function that significantly reduces the computational cost of a single evaluation, compared to the loss function proposed in Ref.~\citen{tonini2025enhanced}, while preserving a relationship between its minimum $(\hat{\boldsymbol \mu}, \hat \Gamma)$ and $(\mathbf u_\mathrm{MAP}, \Gamma_\mathrm{Lap})$ when the forward map $\mathcal F$ in \eqref{eq:model} is affine (\Cref{thm:convergence} and \Cref{thm:convergenceNew}) and extending the result to the non affine case.

We showed the capabilities of the UQ-VAEs to approximate both the mean and covariance of the posterior distribution in the cases of the Poisson problem (\cref{sec:test1}), a non linear problem (\cref{sec:nonLinear}) and a 0D cardiocirculatory model for systemic hypertension (\cref{sec:test2}) and for ventricular septal defect (\cref{sec:test3}).\\
With few observational data, the UQ-VAE estimates of the true parameters are as accurate as the posterior mean. With a large number of observational data, the true parameters were estimated better using the UQ-VAE mean rather than the posterior one, due to the latter being poorly approximated by the sample mean in the case of low regularity of the posterior distribution. A similar behavior was observed when the noise level $\eta$ decreased from high to low. These results show that the UQ-VAE approach is effective in estimating model parameters by solving inverse problems, even for low regular posterior distributions.\\
The main advantage of a trained UQ-VAEs is in performing inverse uncertainty quantification on the estimated parameters in a short computational time. UQ-VAEs estimated posterior covariance matrices, with average relative errors less than $1.6\%$ for the 0D cardiocirculatory model and a test set of $100$ samples. The average relative error on the posterior variance in the case of the Poisson problem for a test set of $100$ samples was less than $11\%$. For this problem, the UQ-VAE estimated about $42,000$ coefficients ($D+D(D+1)/2$, with $D= 17^\mathrm{2}$) for the mean and the covariance matrix, maintaining a good accuracy.

The proposed novel loss function significantly reduced the training time compared to the previous formulation, while also improving accuracy in estimating $(\mathbf u_\mathrm{MAP}, \Gamma_\mathrm{Lap})$. This reduction in computational cost enables the solution of high dimensional Bayesian inverse problems that were previously computationally infeasible.

A limitation of UQ-VAEs arises when the map $\mathcal F$ is approximated by a NN $\psi$. In this case, the decoder is trained on noiseless observational data and the approximation error is also computed using such data, which is only feasible within an in silico framework. In contrast, the encoder is trained using noisy observational data. When the map $\mathcal F$ is known, noiseless observational data are no longer required: an advantage in realistic scenarios where observational data are inherently noisy.\\
Moreover, the decoder is trained using a supervised approach, with paired input and output data, whereas the encoder is trained in a semi-supervised manner, requiring only observational data and prior knowledge about the distribution of the parameters. In real case scenarios, the observational data of a phenomenon are directly measured, while the prior knowledge is represented by a mathematical model describing the phenomenon. The prior mean can be fixed to a reference setting of the model parameters and the prior covariance can be defined by a variation of the prior mean, as we have done in \cref{sec:results}. Therefore, a semi-supervised approach for training the decoder, which requires only observational data, would broaden the applicability of UQ-VAEs, as it would eliminate the need for in silico generated parameters and observational data.\\
Finally, we approximate the posterior distribution using a Gaussian distribution, which is not correct for general maps $\mathcal F$ from parameters to observational data. This limitation can be addressed by using hierarchical VAEs, which enable sampling from more expressive approximate posteriors \cite{kingma2016improved,sonderby2016ladder,vahdat2020nvae}.

\section*{Acknowledgments}
AT, FR, LD and AQ are members of the INdAM group GNCS \say{Gruppo Nazionale per il Calcolo Scientifico} (National Group for Scientific Computing).

AT, FR and LD acknowledge the INdAM GNCS project CUP E53CE53C24001950001.

FR has received support from the project FIS, MUR, Italy 2025-2028, Project code: FIS-2023-02228, CUP: D53C24005440001, \say{SYNERGIZE: Synergizing Numerical Methods and Machine Learning for a new generation of computational models}.

LD acknowledges the support by the FAIR (\say{Future Artificial Intelligence Research}) project, funded by the NextGenerationEU program within the PNRR-PE-AI scheme (M4C2, investment 1.3, line on Artificial Intelligence), Italy.

LD acknowledges the project PRIN2022, MUR, Italy, 2023–2025, 202232A8AN \say{Computational modeling of the heart: from eﬀicient numerical solvers to cardiac digital twins}.

The present research is part of the activities of “Dipartimento di Eccellenza 2023–2027”, MUR, Italy, Dipartimento di Matematica, Politecnico di Milano.

\section*{Competing interests}
The authors declare no competing interests.

\bibliographystyle{plain}
\bibliography{references}

\end{document}